\documentclass[a4paper,11pt]{article}

\usepackage{amsmath,amsfonts,amssymb,mathtools,amsthm}
\usepackage{bbm} 
\usepackage{epsfig,mathrsfs,enumerate,graphics}
\usepackage[usenames,dvipsnames]{color}
\usepackage{dsfont}
\usepackage{graphicx}
\usepackage{wrapfig}
\usepackage{caption}
\usepackage{subcaption}
\usepackage{tikz}
\usepackage{datetime2}
\usepackage[english]{babel}

\usepackage{todonotes}

\usepackage[all]{xy} \entrymodifiers={!!<0pt,0.7ex>+}

\usepackage[hyperindex=true, hidelinks,
					colorlinks=false]{hyperref}

\newcommand{\blue}[1]{\textcolor{black}{#1}}
\newcommand{\black}[1]{\textcolor{black}{#1}}

\newcommand{\ens}{\enspace}

\newtheorem{thm}{Theorem}

\newtheorem{lem}{Lemma}[section]
\newtheorem{prop}[lem]{Proposition}
\newtheorem{cor}[lem]{Corollary}

\theoremstyle{definition}

\newtheorem{rem}[lem]{Remark}
\newtheorem{Def}[lem]{Definition}
\newtheorem{exa}[lem]{Example}

\numberwithin{equation}{section}

\newcounter{parag}[subsection]

\newcounter{parage}

\newcounter{paraga}

\def\al{\alpha}
\def\be{\beta}

\def\Ga{{\Gamma}}
\def\de{\delta}

\def\De{\Delta}
\def\eps{{\varepsilon}}
\def\ka{\kappa}
\def\la{\lambda}
\def\La{\Lambda}

\def\Om{\Omega}

\def\sig{{\sigma}}

\def\th{{\theta}}
\def\Th{\Theta}
\newcommand{\ph}{\varphi}

\newcommand{\ov}{\overline}
\newcommand{\wh}{\widehat}

\newcommand{\adm}{\operatorname{adm}}

\newcommand{\sym}{\operatorname{sym}}
\newcommand{\card}{\operatorname{card}}

\newcommand{\id}{\operatorname{id}}
\newcommand{\ID}{\mathop{\hbox{{\rm Id}}}\nolimits}

\newcommand{\dd}{{\mathrm d}}

\newcommand{\ee}{\mathrm e}

\newcommand{\pa}{\partial}

\newcommand{\ii}{^{-1}}
\newcommand{\ic}{^{\circ(-1)}}
\newcommand{\ti}{\tilde}

\newcommand{\ie}{{\emph{i.e.}}\ }

\newcommand{\eg}{{\it e.g.}\ }

\newcommand{\lhs}{{left-hand side}}
\newcommand{\rhs}{{right-hand side}}

\newcommand{\dst}{\displaystyle}

\newcommand{\C}{\mathbb{C}}

\newcommand{\D}{\mathbb{D}}

\newcommand{\N}{\mathbb{N}}
\newcommand{\Q}{\mathbb{Q}}
\newcommand{\R}{\mathbb{R}}

\newcommand{\Z}{\mathbb{Z}}

\def\cA{\mathcal{A}}
\def\cB{\mathcal{B}}

\def\cF{\mathcal{F}}

\def\cN{\mathcal{N}}

\def\cT{\mathcal{T}}

\DeclarePairedDelimiter\abs{\lvert}{\rvert}%
\DeclarePairedDelimiter\norm{\lVert}{\rVert}%

\def\cS{{\cal S}}

\newcommand{\ord}{\operatorname{ord}}

\newcommand{\htx}[1]{\raisebox{-.09ex}{${\stackrel{
            \raisebox{-.23ex}{$\scriptscriptstyle\wedge$}
          }{#1}
     }$}}

\newcommand{\shtx}[1]{
        \raisebox{0.15ex}[2ex][0ex]{${\stackrel{
            \raisebox{-.12ex}{$\scriptscriptstyle\wedge$}
          }{\scriptstyle#1}
     }$}}

\newcommand\Bbibitem[2]{\bibitem[#1]{#2}}

\newcommand{\defeq}{\coloneqq} 

\newcommand{\col}{\colon\thinspace}          

\newcommand{\gG}{\mathscr G}       
\newcommand{\gP}{\mathscr P}       

\newcommand{\est}{\emptyset}

\newcommand{\begla}{\begin{equation}}
\newcommand{\beglab}[1]{\begin{equation}	\label{#1}}
\newcommand{\edla}{\end{equation}}

\newcommand{\imp}{\ens\Longrightarrow\ens}
\newcommand{\Imp}{\quad\Longrightarrow\quad}

\newcommand{\bul}{\bullet}

\newcommand{\att}{\triangleleft}

\newcommand{\gGnu}{\gG_\nu}
\newcommand{\tgGnu}{\ti\gG_\nu}
\newcommand{\tgGnuid}{\ti\gG_{\nu,\id}}
\newcommand{\gGnuid}{\gG_{\nu,\id}}
\newcommand{\tgPnu}{\ti\gP_{\nu}}

\newcommand{\GLnu}{\operatorname{GL}(\nu,\C)}

\newcommand{\Tree}{\operatorname{Tree}}

\newcommand{\Cznu}{\C[[z_1,\ldots,z_\nu]]}

\newcommand{\dde}{[d]\setminus\De}

\newcommand{\di}[1]{d_{#1}\hspace{-.08em}(i)}

\newcommand{\colb}{\col\hspace{-.25em}}

\newcommand{\ocN}{\wh\cN}

\newcommand{\ppi}{^{(i)}}

\newcommand{\plim}{\lim\limits_{{\longleftarrow}}}
\newcommand{\vph}{\vphantom{\plim\raisebox{.5ex}{$\ocN$}}}

\newcommand{\tmem}[1]{{\em #1\/}}
\newcommand{\tmop}[1]{\ensuremath{\operatorname{#1}}}

\newcommand{\BeunedeuxSetroisequatre}{\begin{tikzpicture}[scale=0.7,baseline=5pt]
  \node (r) at (0,0) {$n_1$};
  \node (a) at (-0.5,1) {$n_2$};
  \node (b) at (0.5,1) {$n_3$};
  \node (c) at (0.5,2) {$n_4$};
    \draw (r) --(a) ;
     \draw (r) --(b) --(c);
 \end{tikzpicture}}
 
  \newcommand{\Becinqesix}{\begin{tikzpicture}[scale=0.7,baseline=5pt]
  \node (r) at (0,0) {$n_5$};
  \node (a) at (0,1) {$n_6$};
    \draw (r) --(a) ;
 \end{tikzpicture}}
 
   \newcommand{\Bedeuxtrois}{\begin{tikzpicture}[scale=0.7,baseline=5pt]
  \node (r) at (0,0) {$n_3$};
  \node (a) at (0,1) {$n_4$};
    \draw (r) --(a) ;
 \end{tikzpicture}}

\newcommand{\BeBeunedeuxSetroisequatreSecinqesix}{\begin{tikzpicture}[scale=0.7,baseline=5pt]
  \node (r) at (0,0) {$n$};
 \node (a) at (-0.5,1) {$n_1$};
 \node (b) at (-1,2) {$n_2$};
  \node (c) at (0,2) {$n_3$};
  \node (d) at (0,3) {$n_4$};
   \node (e) at (0.8,1) {$n_5$};
    \node (f) at (0.8,2) {$n_6$};
    \draw (r) --(a)--(b);
     \draw (a) --(c)--(d);
     \draw (r) --(e) --(f);
\end{tikzpicture}}

\newcommand{\datestampa}{{\small{File:\ens\hbox{\tt\jobname.tex}
\ens \today}}}
\newcommand{\datestamp}{\small{File:\ens\hbox{\tt\jobname.tex}
\ens \DTMnow}}

\newlength{\locwidth}
\newlength{\locheight}
\newcommand{\ppb}[1]{\settowidth{\locwidth}{#1}%
  \settoheight{\locheight}{#1}%
\raisebox{.3\locheight}{$\left(\parbox[c]{\locwidth}{#1}\right)$}}

\newenvironment{itquote}
  {\begin{quote}\itshape}
    {\end{quote}\ignorespacesafterend}

\providecommand\given{}  
\newcommand\SetSymbol[1][]{\nonscript\:#1\vert\nonscript\:\allowbreak}
\DeclarePairedDelimiterX\Set[1]\lbrace\rbrace{%
  \renewcommand\given{\SetSymbol[\delimsize]}#1}

\begin{document}

\thispagestyle{empty}
\begin{center}
{\bf \huge
Explicit linearization of multi-dimensional \\[1ex] germs and vector fields through\\[1ex] Ecalle's tree expansions}

\vspace{1.5cm}

Fr\'ed\'eric Fauvet, Fr\'ed\'eric Menous, David Sauzin

\vspace{1.5cm}

\today\footnote{\datestamp}

\end{center}
\vspace{2cm}

\begin{abstract}
  \black{We revisit the Siegel linearization problem by means of
    explicit formulas of non-recursive type for the linearizing
    transformations of a non-resonant holomorphic germ of diffeomorphism
    at a fixed point, or a non-resonant holomorphic germ of vector field
    at a singular point, in any complex dimension. These formal
    expressions, that have the same shape for dynamical systems with
    discrete or continuous time and were obtained long ago by
    J.~Ecalle through his ``arborification'' technique, are derived
    here by a more direct application of his tree-based combinatorics
    called ``armould calculus".
    They involve a certain family of differential operators indexed by decorated
    trees and forests. 
    A combinatorial study of these differential
    operators allows us to recover, under Bruno's arithmetical
    condition, the best known estimates for the domains of convergence
    of the holomorphic linearizing changes of variables in terms of the
    value of the Bruno series, including a new precise dependence with
    respect to the dimension of the problem.}
    %
%
\end{abstract}

\newpage

\setcounter{tocdepth}{2}
\tableofcontents

\vspace{1cm}

\begin{center} \rule{3cm}{.1pt} \end{center}

\vspace{1cm}








\section{Introduction}\label{intro}
In the present article, we revisit the classical \black{Siegel
  linearization problem: the conjugacy to their linear parts} of
complex analytic dynamical systems at a singular point of
$\mathbbm{C}^{\nu}$ (any $\nu \geq 1$), with discrete or
continuous time. By using a part of Ecalle's ``mould calculus'' 
involving tree-based combinatorics, \black{one can} 
give explicit formulas for the normalizing transformations, which are
of the same type for diffeomorphisms and vector fields; \black{we will
  use them to derive
estimates} for their domains of convergence.
%

  \subsection{\blue{Overview}}

It has been known since Poincar{\'e} that any local analytic 
diffeomorphism~$g$ that fixes the origin of $\mathbbm{C}^{\nu}$ and
has a linear part of the form 
\beglab{eqdefRq}
\dd g(0) = R_q \col (z_1,\ldots, z_\nu)
\mapsto (q_1 z_1,\ldots, q_\nu z_\nu),
  \edla
  with a spectrum $q = (q_1,\ldots,q_\nu) \in (\C^*)^\nu$ satisfying a
  so-called \emph{non-resonance} condition (see Section~\ref{seclingermnu}),
  is formally conjugate to $\dd g (0)$, with a unique
  tangent-to-identity linearizing 
  transformation~$h$. Formulas of a recursive type for~$h$
%
have been used for more than a century and~$h$ has been
proved convergent when~$q$ is in the so-called Poincar{\'e} domain,
that is the eigenvalues $q_1, \ldots, q_{\nu}$ 
are either all in the open unit disc $\D$ of~$\C$, 
or all outside $\overline{\D}$.
%
%
If the
spectrum~$q$ is non-resonant but not in the Poincar\'e domain, then
some \black{nonzero} integer combinations of the
$q_i$'s, appearing as denominators in the recursive formulas, can be
arbitrarily close to
$0$ and may thus prevent the convergence
of~$h$: this is the ``small denominator problem''.

It is only in 1942 that C. Siegel, in the celebrated paper \cite{Siegel}, was able to
control these quantities and reach a result of analytic linearization
in dimension~$1$
%
%
under a suitable diophantine hypothesis on the 
eigenvalue~$q$. He later generalized this result to
higher-dimensional germs, with similar techniques under an arithmetical
condition of the same type on the spectrum~$q$. 
In the real $\mathcal{C}^{\infty}$ category, Sternberg
\cite{Stern} used another strategy to prove that a smooth non-resonant local
diffeomorphism fixing the origin of $\mathbbm{R}^{\nu}$ is
linearizable, with smooth changes of coordinates.  Siegel's results
were sharpened by Bruno in 1971 and R\"ussmann in 1977, with the
introduction of weaker arithmetical conditions on the spectrum $q$,
expressed by the convergence of a series $\mathcal{B} (q)$ of
non-negative real numbers (see \cite{Bruno2,russmann}, the survey \cite{Herman} and
Section~\ref{seclingermnu}).

Parallel results \black{have been obtained} 
for local analytic vector fields vanishing at the origin.
%
%
For an analytic vector field~$V$ whose $1$-jet is
\begla
V_{\lambda}\defeq \sum_{i=1}^\nu \lambda_i z_i
\frac{\partial}{\partial z_i},
\edla
with non-resonant spectrum
$\lambda=(\lambda_1,\ldots,\lambda_\nu)\in\C^\nu$,
Bruno had by 1971 (\cite{Bruno}) introduced the arithmetical condition
on~$\lambda$ that now bears his name, involving a series
$\cB_{v.f.} (\la)$:
if $\cB_{v.f.} (\la) < \infty$, then the vector field is analytically
linearizable (see Section~\ref{seclingermnu} for precise definitions
and statements).

\color{black}
As for quantitative estimates of the domains of convergence of the linearizing
transformation~$h$ in terms of the Bruno series,
a landmark was Yoccoz's celebrated 1987 work published in \cite{Yoccoz},
giving a lower bound of the order of $e^{-\cB(q)}$ for
one-dimensional holomorphic germs (instead of the
previously known $e^{-2\cB(q)}$), 
and establishing the optimality of Bruno's condition in that context
by means of his novel geometric renormalization theory.
More recently, the articles \cite{GLS} and \cite{GM}, dealing with
diffeomorphisms and vector fields respectively in arbitrary dimention,
brought improvements in the estimates of the domains of convergence in
line with Yoccoz's lower bound.
%
%
The results that we present in this paper compare well with these.
%
\color{black}

In the present work, we essentially follow a scheme \black{designed by}
Ecalle (see \cite{EcNonAbord} for instance) 
to tackle very general classes
of problems,
that provides a remarkable way of contructing
the unique tangent-to-identity linearizing transformation~$h$:

\noindent
-- the object we focus on is the substitution operator $C_h : \varphi \mapsto
\varphi \circ h$ associated to~$h$;\\[.8ex]
-- the operator $C_h$ is {\tmem{a priori}} expressed as an expansion over a
set of $\mathbbm{Z}^{\nu}$-decorated forests~$F$:
\beglab{eqmldexpaCh}
C_h = \sum_F S^F D_F,
  \edla
  where the $S^F$ are complex numbers that depend only on the spectrum and the
  $D_F$ are 
  differential operators in~$\nu$ variables that depend only
on the nonlinear part of the data.
The collection of numbers $(S^F)$ is called an
armould\footnote{Armoulds were formerly called arborescent moulds or
  arbomoulds by Ecalle. Coarmoulds were formerly
  called arborescent comoulds or coarbomoulds.}
and we denote it by $S^{\bullet}$;
the collection of operators $(D_F)$ is called a coarmould and we denote it by $D_{\bullet}$.
We will obtain explicit formulas for~$S^{\bullet}$ and~$D_{\bullet}$ that
are essentially the same for diffeomorphisms and for vector fields
({\tmem{mutatis mutandis}} as regards the armould~$S^{\bullet}$), and
hence explicit formulas for~$h$ itself,
from which we will derive analyticity results
under the Bruno condition
almost at one stroke.

More precisely, denoting by~$\D_b^\nu$ the polydisc of
radius~$b$ centred at~$0$ in~$\C^\nu$,
if we take a diffeomorphism
%
%
or a vector field 
of the form
\beglab{eqdefgorV}
g = R_q \circ (\tmop{id} +a)
\quad \text{or}\quad
V = V_{\lambda} + a,
\edla
where the nonlinear part~$a$ is analytic \black{and bounded}
on~$\D_b^\nu$, 
we obtain (in Theorem~\ref{thm:bounds} of Section~\ref{secthmbounds}) the following explicit
lower bound for the radius of a polydisc in~$\C^\nu$ in which the corresponding
linearizing transformations~$h$ are convergent:
\begla
b \, \frac{e^{- \gamma -\mathcal{B}}}{\nu (4 M \nu + 2)},
\edla
where~$M$ is the supremum of $\frac{1}{b} a$ on
$\D_b^\nu$, 
$\gamma$ is some explicit numerical constant and $\mathcal{B}$ stands for
$\mathcal{B} (q)$ or $\cB_{v.f.} (\la)$. 
This is a refinement of the results of \cite{GLS} and \cite{GM}, as
explained in Section~\ref{ssec:estim}.

We refer to \cite{FM,FMS18} for explanations on the relevance of a
combinatorial system
%
involving trees and forests in these questions and related algebraic
constructions, yet all the necessary material for our purposes will be
presented below.\footnote{%
  Trees were introduced by
  Cayley \cite{Cayley} in differential calculus
  and, in 1988, by Eliasson 
for KAM theory (\cite{Eliasson}).
Another tree
formalism was used in \cite{Carletti} for the multidimensional Siegel problem.%
}
In \cite{FMS18}, we had implemented Ecalle's armould
formalism for the linearization of non-resonant diffeomorphisms in one
complex dimension and shown how it yields explicit formulas instead of
the classical recursive ones, that enabled us to improve upon the
existing results in the literature;
%
going from one-dimensional to multidimensional brings a number of technical
complications already at the formal level.
The final section of \cite{FM} touches upon the analytic aspects of the multidimensional
linearization problem, but the focus of that paper is definitely algebraic.
As explained in \cite{FM},
Ecalle's theory of arborification (\cite{EcNonAbord}) can be formulated within a framework
of Hopf algebras,
 involving
 notably the Connes-Kreimer Hopf algebra (\cite{CK}) on the armould side and the
Grossman--Larson one for coarmoulds. This rich underlying formalism
is however not needed in the present article, the core of which resides in the
technical algebraic and combinatorial constructions of Section~\ref{sec:vanish} together with the subsequent
sharp estimates derived in Section~\ref{sec:estim}.

\subsection{\black{Plan}}

The paper is organized as follows:
\\[.8ex]
\noindent --
In Section~\ref{seclingermnu}, we set the stage and introduce our
notations for the (formal or holomorphic) linearization problems for
diffeomorphisms or vector fields, for which we briefly review the
literature, and for our Bruno series, which we compare with the ones used
in \cite{GLS} and \cite{GM}.
\\[.8ex]
\noindent --
In Section~\ref{sec:formal}, we introduce terminology pertaining
to Ecalle's armould calculus, define and prove some key properties of
the homogeneous operators~$D_F$ of~\eqref{eqmldexpaCh}, and finally
prove---up to Lemma~\ref{propnullDcFplus}, whose proof is deferred to
Section~\ref{sec:vanish}---that formula~\eqref{eqmldexpaCh} makes
sense and gives the solution to the formal linearization problem when
using the appropriate armould~$(S^F)$:
see Theorems~\ref{thmformaldiffeonu} and~\ref{thmformalfieldnu}
for diffeomorphisms and vector fields of the form~\eqref{eqdefgorV},
with remarkably compact formulas for the relevant coefficients~$S^F$.
The definition of both~$D_F$ and~$S^F$ is due to Ecalle \cite{EcNonAbord}.
\\[.8ex]
\noindent --
Section~\ref{sec:analytic} shows how to deduce from the explicit
tree expansions~\eqref{eqmldexpaCh} analyticity results for the
linearizing transformations~$h$ under the Bruno condition;
the upshot is Theorem~\ref{thm:holomorphic}, equivalent to Bruno's and
R\"ussmann's results, and its refinement, Theorem~\ref{thm:bounds}, which gives
precise estimates for the domains of convergence.
Proofs are carried out in a self-contained manner up to two results
whose proofs are deferred to later sections, namely
an algebraic cancellation phenomenon explained in
Section~\ref{sec:vanish}
and the armould estimates given in Section~\ref{sec:estim}.
\\[.8ex]
\noindent --
Section~\ref{sec:vanish} is of algebraic nature, 
with emphasis on the dichotomy between those forests~$F$ for
which~$D_F$ is always the zero operator and the others.
For the latter ones, the so-called ``admissible cuts'' obey
Theorem~\ref{th:cutvanish}, that we state and prove,
and from which we deduce the instrumental
Proposition~\ref{prop:annul2} (required in the proof of
Lemma~\ref{propnullDcFplus} and some passages of
Sections~\ref{sec:analytic} and~\ref{sec:estim}).
\\[.8ex]
\noindent --
Section~\ref{sec:estim} contains the crucial estimates for the relevant
armoulds $S^{\bullet}$, in the form of Theorem~\ref{thm:estimvf} for
vector fields and Theorem~\ref{thm:estimq} for diffeomorphisms,
whose proofs only require a little combinatorial effort and
manipulations in which the Bruno series appears in the most natural way.
\\[.8ex]
\noindent --
Appendix~\ref{secappcomputDF} gives closed-forms expressions for the
components of the coarmould $D_{\bullet}$.

\medskip

We hope to convince the reader that the appropriate combinatorial
framework allows one to reduce the analytic estimates to the bare
minimum.
Note that Theorems~\ref{thmformaldiffeonu} and~\ref{thmformalfieldnu}
were proved by Ecalle in \cite{EcNonAbord} by means of the so-called
``arborification'' technique, using a preliminary mould expansion
prior to reaching the tree expansion~\eqref{eqmldexpaCh}, whereas
our method is more direct.
Theorem~\ref{thm:holomorphic} is also given in \cite{EcNonAbord}, but
the proof there is very concise and we do not know how to obtain the
necessary armould estimates without making use of our
Proposition~\ref{prop:annul2};
the situation is somewhat similar with \cite{FM}, which mainly
  focuses on the Hopf-algebraic aspects of Ecalle's arborification and
  does not go into all the analytic details for the linearization problem.

  \subsection{\black{Outlook}}

\black{Beyond the obtention of explicit formulas for the 
linearization of non-resonant local dynamical systems
and estimates for the domain of analyticity under Bruno's arithmetical
condition, the present paper paves the way for treating interesting further
questions.
First, since Bruno's condition is not believed to be optimal for the Siegel
linearization problem when the dimension~$\nu$ is larger than~$1$, one
could try to replace it by a weaker condition, in the spirit of
the recent preprint \cite{ArCor}.
One could also try to adapt the formalism to tackle the more
complicated KAM problem, thus enhancing the semi-direct method of
\cite{Eliasson}, as is claimed to be feasible \eg in \cite{EcVallet}.}

Another avenue of research would be to implement effectively the research
program outlined in \cite{EcNonAbord} for a resonant spectrum,
starting with the simply resonant case (\ie
Definitions~\ref{defqnonres} or~\ref{deflanonres} are violated but all
resonance relations can be deduced from one of them).
Then, generically, the normal form is not reduced to the linear
part and divergence of resurgent type occurs in the formal normalizing transformations
(see \cite{mouldSN} for
two-dimensional saddle-node vector fields, and \cite{costin},
\cite{Kamimoto} for partial results in higher dimension); in general, the coexistence of small
denominators and resonance requires an arithmetical condition for a
full analytic treatment of the normalization problem.

\black{Even in the purely resonant case, in the absence of small denominators,
  there is a highly non-trivial \emph{analytic classification problem} giving
  rise to functional moduli. The armould apparatus seems to be an
  efficient constructive tool to tackle the ``synthesis'' question,
  which asks to construct, for an arbitrary element of the moduli space, a
  local analytic dynamical system belonging to the analytic
  equivalence class associated with that element.
  In the context of saddle-node singularities (Ecalle-Martinet-Ramis moduli spaces),
  convincing arguments have already been developed in
  \cite{EcTwisted}
(see \cite{FFparalog} for a review):
  the idea is to obtain the sought-after germ of analytic system in
  the form of a convergent tree expansion involving the so-called paralogarithmic monomials.
More generally, a problem like Kontsevich-Soibelmann's resurgence
conjecture \cite{KS22} should fall within the scope of the armould method.}

Still another possible application of the armould formalism would be in
line with P\"oschel's paper \cite{poschel},
where the author shows how to derive Sternberg's as well as other results
with data involving various classes of ultradifferentiable series (in
particular the Gevrey classes) from the recursive formulas used in the
analytic class, involving new insights on classes that are stable by
composition;
within our formalism, this could probably achieved in a more direct
way,
by means of explicit tree expansions instead of recursive formulas,
treating on an equal footing 
diffeomorphisms and vector fields.

We postpone such developments for future research.
Note that one benefit of working at the level of substitution operators
expressed as sums of elementary components,
as Ecalle's armould calculus invites us to do,
is that, beyond spaces of formal
series with growth conditions on their coefficients, these techniques very
naturally lend themselves to dealing with data belonging to some given
{\tmem{functional spaces}}, with estimates expressed through relevant
semi-norms, \eg the ones that can be introduced for
parameter-dependent problems.

\section{The linearization problem for germs of holomorphic diffeomorphisms or vector fields in higher dimension}
\label{seclingermnu}


Throughout this paper we use the notations
\begla
  \N \defeq \{0,1,2,\ldots\}, \quad
  \N^* \defeq \{1,2,\ldots\}
  \quad \text{and}\quad
  [\nu] \defeq \{1,\ldots,\nu\}
  \ens\text{for}\ens \nu\in\N^*.
  \edla

\subsection{The linearization problems in the groups $\protect\gGnu$
  and $\protect\tgGnu$}  \label{seclinpbgr}

We define~$\gGnu$ as the group (for composition) of germs of holomorphic
diffeomorphisms of $(\C^\nu,0)$,
and~$\tgGnu$ as its formal counterpart: the latter is formed of all
$\nu$-tuples $g = (g_1,\ldots,g_\nu)$ of formal series 
$g_1(z_1,\ldots,z_\nu), \ldots, g_\nu(z_1,\ldots,z_\nu) \in \Cznu$ without constant
term
such that the $1$-jet of~$g$ (which is thus a linear endomorphism
of~$\C^\nu$) is invertible.
One can view~$\gGnu$ as a subgroup of~$\tgGnu$, the group law being
induced by the composition of formal series without constant term. 
Inversion in~$\tgGnu$ will be denoted by $g \mapsto g\ic$.

Given $g \in \tgGnu$, its $1$-jet is also called \emph{linear part
  of~$g$} and denoted by $\dd g(0)$, so that
\beglab{eqgddgzerof}
g = \dd g(0) \circ f, \qquad f \in \tgGnuid,
\edla
where $\tgGnuid \subset \tgGnu$ is the subgroup of tangent-to-identity
formal diffeomorphisms.
The general element of $\tgGnuid$ can be written
\beglab{eqgenftgGnuid}
f=(f_1,\ldots,f_\nu), \qquad
f_i(z_1,\ldots,z_\nu) = z_i + \sum_{m\in\N^\nu,\ \abs m\ge2} \ti a_{i,m} z^m
\quad\text{for $i \in [\nu]$.}
\edla
%
%
Here, $(\ti a_{i,m})$ is any collection of complex numbers and we use
the standard notations for multi-integers:
\begla
\abs{m} \defeq m_1 + \cdots + m_\nu \quad\text{and}\quad
z^m \defeq z_1^{m_1} \cdots z_\nu^{m_\nu}
\ens\text{for}\ens m = (m_1,\ldots,m_\nu) \in \Z^\nu.
\edla
We can rewrite~\eqref{eqgenftgGnuid} as
\begla
f = \tmop{id} + a, \quad a \in\tgPnu
\edla
by introducing the set~$\tgPnu$ of all
$\nu$-tuples $a=(a_1,\dots,a_\nu)$ of power series such that, for $i$
in $[\nu]$,
$\displaystyle a_i(z)=\sum_{m\in\N^\nu} \ti a_{i,m} z^m \in \Cznu$ has
$\ti a_{i,m}=0$ when $\abs{m}\le 1$.
%
%
We will refer to $a=f-\tmop{id}$ as to the nonlinear part of any
$g\in\tgGnu$ of the form~\eqref{eqgddgzerof}.


Given $g\in\tgGnu$, the \emph{formal linearization problem} asks for an
$h\in\tgGnu$ conjugating~$g$ and a linear map $L\in \GLnu$.
If such an~$h$ exists, then $\dd h(0)$ establishes a conjugacy
between $\dd g(0)$ and~$L$, hence $h\circ \dd h(0)\ic \in \tgGnuid$ is a conjugacy
between~$g$ and $\dd g(0)$;
therefore, the formal linearization problem can be reduced to
searching for $h \in \tgGnuid$ such that $g \circ h = h \circ \dd g(0)$.

The \emph{holomorphic linearization problem} asks for a holomorphic
conjugacy $h\in\gGnu$ between~$g$ and a linear map under the
assumption that $g\in\gGnu$.


\begin{Def}   \label{defqnonres}
We say that $q = (q_1,\ldots,q_\nu) \in (\C^*)^\nu$ is
\emph{non-resonant} if
\[
i \in [\nu], \;\; m \in \N^\nu \;\;\text{and}\;\; \abs{m}\ge 2
\Imp q^m-q_i \neq 0,
\]
where $q^m = q_1^{m_1} \cdots q_\nu^{m_\nu}$.
We say that a formal diffeomorphism~$g$ is \emph{non-resonant} if
%
%
its linear part $\dd g(0)$ is diagonalizable with eigenvalues
$q_1,\ldots,q_\nu$ giving rise to a non-resonant $\nu$-tuple~$q$.
\end{Def}


%
Performing a preliminary linear change of
coordinates, we can assume without loss of generality that $\dd g(0)$ is
the diagonal linear map $R_q \in \GLnu$ defined by~\eqref{eqdefRq}.
%
%
Therefore, we will assume $g = R_q \circ f$ with $q$ non-resonant and $f \in \tgGnuid$.
It is well-known that, in that case, the formal linearization problem admits a unique
solution, \ie there is a unique 
\beglab{eqdefhim}
h = (h_1,\ldots,h_m) \in \tgGnuid, \quad
h_i(z) = z_i + \sum_{m\in\N^\nu,\ \abs m \ge2} h_{i,m} z^m
\edla
such that 
\begla \label{eqqconj}
\boxed{\ens \vph
  R_q \circ f \circ h = h \circ R_q.
  \ens}
\edla
\'Ecalle's tree formalism will provide an explicit formula for
this conjugating formal diffeomorphism~$h$, see
Theorem~\ref{thmformaldiffeonu} of Section~\ref{sec:formal}.


It is a classical result by R\"ussmann \cite{russmann} and Bruno \cite{Bruno2} that the holomorphic
linearization problem admits a solution under the following
reinforcement of the non-resonance condition.
For a given non-resonant $q = (q_1,\ldots,q_\nu) \in (\C^*)^\nu$,
let $\la=(\la_1,\dots\la_{\nu})\in \C^\nu$ be such that
$q_i=e^{2i \pi \la_i}$ for every $i\in [\nu]$.
We use the notation
\[
  m\cdot \lambda\defeq m_1 \lambda_1 +\dots +m_\nu \lambda_\nu
  \quad \text{for $m\in \N^\nu$}
\]
and notice that the non-resonance condition says that
$m \cdot \la -\la_i \notin\Z$ when $\abs{m}\ge2$.
We set
\beglab{eqRuss} 
%
%
\Omega(k) 
%
%
%
\defeq \min \{1\}\cup \Set[\big]{ d(m \cdot \la -\la_i,\Z) \given i \in [\nu],
\  m \in \N^\nu \ \text{such that}\  2 \leq \abs{m}\leq k+1 }
%
\edla 
for each $k\in\N^*$. This is a positive number independent on the choice of $(\la_1,\dots\la_{\nu})$.

\begin{Def}\label{def:Brunodiff}
We say that a non-resonant $q = (q_1,\ldots,q_\nu) \in (\C^*)^\nu$ satisfies the
Bruno condition if $\cB(q) < \infty$, where 
\begin{equation} \label{eqBrunoseries}
\cB(q) \defeq \sum_{k=1}^{\infty} \left( \frac{1}{k}-\frac{1}{k+1}\right) \log\frac{1}{\Omega(k)}
\qquad \text{(series of non-negative terms).}
\end{equation}
We say that a formal diffeomorphism~$g$ satisfies the Bruno condition
%
%
if its linear part $\dd g(0)$ is diagonalizable with eigenvalues
$q_1,\ldots,q_\nu$ giving rise to a non-resonant $\nu$-tuple~$q$
satisfying the Bruno condition.
\end{Def}

We will explain how \'Ecalle's tree formalism allows one to recover
Bruno-R\"ussmann's result---see Theorem~\ref{thm:holomorphic}(i) of Section~\ref{sec:analytic}.
%
%
Similar to the one-dimensional case (\cite{FMS18}), the method consists in 
\begin{itemize}
\item
  solving explicitly the formal linearization problem by writing the
  composition operator~$C_h$ of the formal linearization in the form
  of a series of operators
  associated with~$f$ and indexed by forests (namely a tree expansion), 
\item
  and, in the case of the holomorphic linearization problem, providing
  suitable bounds by exploiting the Bruno condition and the holomorphy of~$f$.
\end{itemize}
The same techniques allow one to prove similar results
in the case of holomorphic vector
fields.

\subsection{The linearization problem for germs of holomorphic vector
  fields}  \label{seclinpbvf}

Near~$0$ as an equilibrium point, we consider formal or holomorphic
vector fields with diagonal linear part,
%
%
\ie $\nu$-tuples $V = (V_1,\ldots,V_\nu)$ of formal series of the form
%
%
\beglab{eqcoefftiaiminVi}
V_i(z_1,\ldots,z_\nu)=\la_i z_i +\sum_{m \in \N^\nu,\ |m|\geq 2} \ti a_{i,m} z^m
  \in \Cznu
  \quad \text{for} \ens i=1,\ldots,\nu.
  \edla

Such a vector field decomposes in a linear part $V_\la (z)=(\la_1
z_1,\dots ,\la_\nu z_\nu)$ and a nonlinear one $a=V-V_\la\in \tgPnu$
and, as in the case of diffeomorphisms, we have a linearization problem that can be easily expressed with differential systems associated to such vector fields.

Consider the formal system of differential equations:
\[
\frac{dz_i}{dt}=V_i(z)=\la_i z_i+a_i(z) \quad \text{for} \ens i=1,\ldots,\nu.
\]
The formal linearization problem asks for an $h \in \tgGnuid$ such that, under the formal change of coordinates $z=h(\tilde{z})$ the system becomes
\[
\frac{d \tilde{z}_i}{dt}=\la_i \tilde{z}_i \quad \text{for} \ens i=1,\ldots,\nu,
\]
\ie $h$ is a formal identity-tangent diffeomorphism that formally
conjugates $V$ to its linear part.
The conjugacy equation is obtained by writing
%
\begin{equation}\label{eq:vecConj}
  \frac{d z_i}{dt}=\sum_{j\in [\nu]}
  \frac{d\tilde{z}_j}{dt}\frac{\partial h_i}{\partial
    \tilde{z}_j}(\tilde{z})=
  \boxed{\vph \sum_{j\in [\nu]} \la_j\tilde{z}_j\frac{\partial h_i}{\partial \tilde{z}_j
  }(\tilde{z})=V_i(h(\tilde{z}))}
  =V_i(z)
  \quad \text{for} \ens i=1,\ldots,\nu.
\end{equation}
As for holomorphic germs of diffeomorphisms, there is a holomorphic linearization problem that asks to find a holomorphic conjugacy $h\in  \gGnu$ between a holomorphic vector field and its linear part.

\begin{Def}  \label{deflanonres}
The $\nu$-tuple $\la = (\la_1,\ldots,\la_\nu)$ is the \textit{spectrum} of $V=V_\la +a$ and we say that $\la$ is
\emph{non-resonant} if
\[
i \in [\nu], \;\; m \in \N^\nu \;\;\text{and}\;\; \abs{m}\ge 2
\Imp m\cdot \la-\la_i \neq 0.
\]
%
%
We say that a formal vector field~$V$ is \emph{non-resonant} if its spectrum is non-resonant.
\end{Def}

It is well-known that \black{any non-resonant formal vector field can
  be formally linearized: there is a unique $h\in\tgGnu$
  solving~\eqref{eq:vecConj}}---see Theorem~\ref{thmformalfieldnu} of
Section~\ref{sec:formal} for a proof using Ecalle's techniques.

From the holomorphic point of view, the vector field analog of
Bruno-R\"ussmann's result is 
Bruno's theorem \cite{Bruno},
\black{which says that an analytic vector field can be analytically
  linearized under the folowing reinforcement of the non-resonance condition.}
For a given non-resonant spectrum $\la=(\la_1,\ldots,\la_\nu)$ and any positive integer $k$, let
\begin{equation}\label{eqBru}
\Omega_{v.f.}(k)\defeq  \min \{1\}\cup \Set[\big]{ |m\cdot \la -\la_i| \given i
\in [\nu], \  m \in \N^\nu \ \text{such that}\  2 \leq \abs{m} \leq k+1 }.
\end{equation}


\begin{Def} \label{def:BrunoVF}
We say that a non-resonant $\la=(\la_1,\ldots,\la_\nu) \in \C^\nu$ satisfies the
Bruno condition if $\cB_{v.f.}(\la) < \infty$, where 
\[
\cB_{v.f.}(\la) \defeq \sum_{k=1}^{\infty} \left( \frac{1}{k}-\frac{1}{k+1}\right) \log\frac{1}{\Omega_{v.f.}(k)}
\qquad \text{(series of non-negative terms).}
\]
We say that a formal vector field $V$ satisfies the Bruno condition
if its spectrum satisfies the Bruno condition.
\end{Def}

%

As in the case of holomorphic germs, Ecalle's tree formalism allows
one to get explicit formulas and to recover Bruno's theorem, see
Theorem~\ref{thm:holomorphic}(ii) of Section~\ref{sec:analytic}.

\subsection{Estimates and comparisons}\label{ssec:estim} 


As mentioned in the introduction, Theorem~\ref{thm:bounds} of
Section~\ref{secthmbounds} below shows how Ecalle's formalism allows one to get
estimates on the radius of convergence of the conjugacy map that
%
%
recover the best known estimates due to \black{\cite{Yoccoz},} \cite{GLS} and \cite{GM}, with
a new explicit dependence with respect to the dimension.

In both linearization problems, starting with a holomorphic
diffeomorphism or a vector field of the form~\eqref{eqdefgorV} with
nonlinear part~$a$ analytic \black{and bounded} on the polydisc
$\D_b^\nu$ of radius $b$,
we obtain the following lower bound for the radius of a polydisc in
which the conjugacy map~$h$ is convergent:
\[
  b \, \frac{e^{- \gamma -\mathcal{B}}}{\nu (4 M \nu + 2)}
  \quad \text{with}\ens \mathcal{B}=\mathcal{B}(q)
  \ens \text{or} \ens\cB_{v.f.}(\la),
  \]
where $M$ is the supremum of $\frac{1}{b}\abs a$ on $\D_b^\nu$ 
and $\gamma$ is some explicit numerical constant.
%

This is in accordance with the results of \black{\cite{Yoccoz},} \cite{GLS} and \cite{GM},
which were the first articles to give a lower bound of the form
$e^{-K-\mathcal{B}}$ with~$K$ depending on the nonlinear part~$a$ but
not on~$\la$
(prior to them only lower bounds of the form $e^{-K-2\mathcal{B}}$ were known).
Note that, for vector fields, our Bruno series $\cB_{v.f.}(\la)$ is the same as in \cite{GM},
%
but for diffeomorphisms 
there is a slight discrepancy between our Bruno series $\cB(q)$ and
\black{the Bruno function used in \cite{Yoccoz} when $\nu=1$, for which
  case we refer to Section~\ref{seconedim}, or with
  the series used in \cite{GLS} when $\nu\ge2$.}
The latter reference makes use of the same series as in~(\ref{eqBrunoseries}) except for the
function measuring the size of the small denominators.
We measure this size by the function~$\Omega$ defined in
equation~(\ref{eqRuss}), as in \cite{Bruno2}, whereas
in \cite{GLS} $\Omega(k)$ is replaced by
\begin{equation}\label{eqGls}
  \al(k)\defeq \min \{1\}\cup \Set[\big]{ | e^{m \cdot \la -\la_j}-1|
  \given j \in [\nu], \  m \in \N^\nu \ \text{such that}\  2 \leq \abs{m}\leq k+1 }.
\end{equation}
%
However, these functions are of the same order 
in the sense that there exist universal constants $c,C>0$ such that
\[
c\al(k)\leq \Omega(k) \leq C \al(k) \quad \text{for all}\; k>0.
\]
Therefore, the two Bruno series 
\[
\cB(q) = \sum_{k=1}^{\infty} \left( \frac{1}{k}-\frac{1}{k+1}\right)
\log\frac{1}{\Omega(k)}
\quad\text{and} \quad
\sum_{k=1}^{\infty} \left( \frac{1}{k}-\frac{1}{k+1}\right) \log\frac{1}{\al(k)}
\]
have the same nature and, in case of convergence, their difference is
bounded by a constant depending only on~$c$ and~$C$.

For the sake of completeness, let us mention that a third function
$\varepsilon$ appears in \cite{russmann}, \cite{poschel} and
\cite{raissy}, defined by 
\begin{equation}\label{eqrais}
  \varepsilon(k)\defeq \min \{1\}\cup \Set[\big]{ | q^m-q_j| \given j \in
  [\nu], \  m \in \N^\nu \ \text{such that}\  2 \leq \abs{m}\leq k+1 }. 
\end{equation}
This function is of the same order as the previous ones, so that the
corresponding Bruno condition is equivalent. Note also that there are
different types of Bruno series that are proved to be equivalent in
the appendix of \cite{raissy}.

\subsection{Remark on the one-dimensional case}   \label{seconedim}


For a non-resonant $q\in (\C^*)^\nu$, in any dimension~$\nu$,
since~\eqref{eqRuss} defines a monotonic non-increasing sequence of
positive real numbers, there is a unique increasing sequence of
positive integers $(k_j)_{j\ge0}$ such that $k_0=1$ and
\[
  k_j \le k < k_{j+1} \imp
  \Om(k_j) = \Om(k) > \Om(k_{j+1})
  \quad
  \text{for all $j\ge0$.}
  \]
  In general, the above sequence $(k_j)$ may be finite. If~$q$ is in the Siegel domain, then $\Om(k)$ tends to~$0$ as
  $k\to\infty$, hence the sequence is infinite and the series~$\cB(q)$
  defined in~\eqref{eqBrunoseries} can be rewritten
  \beglab{eqaltcBq}
  \cB(q) =\sum_{j=1}^{\infty} \left( \frac{1}{k_j}-\frac{1}{k_{j+1}}\right) \log\frac{1}{\Omega(k_j)}.
  \edla

  \begin{lem}
    Let $\nu=1$ and $q=e^{2i \pi \la}$ with $\la\in\R\setminus\Q$
    (thus $q$ is in the Siegel domain).
    The above sequence $(k_j)$ then coincides with the sequence
    $(q_j)$ of the denominators of the convergents of the continued
    fraction expansion of~$\la$,
    and
    \begla
    \Om(q_j) = d(q_j\la,\Z)
    \quad \text{for all $j\ge0$.}
    \edla
  \end{lem}

  \begin{proof}
    Let us denote by $(p_j/q_j)$ the sequence of convergents of~$\la$.
    Their well-known properties, in particular the law of the best approximation, imply that
    $d(n\la,\Z) \ge d(q_j\la,\Z)$ for $1\le n<q_{j+1}$.
    Since $\Om(k) =  \min \Set[\big]{ d(n \la,\Z) \given 1 \leq n\le k
    }$, we get
    \[
      1\le k<q_{j+1} \imp \Om(k) =\Om(q_j)=d(q_j\la,\Z),
    \]
    which is equivalent to the desired result.    
  \end{proof}

  \begin{cor}
    The difference between $\cB(q)$ and the Bruno function
    \beglab{eqdefBrfcn}
    \sum_{j=0}^{\infty} \frac{1}{q_j} \log q_{j+1}
    \edla
    is bounded by a universal constant.
  \end{cor}

  \begin{proof}
    Another well-known property of the convergents of~$\la$ is
    \[
      \frac 1 {2q_{j+1}} < \frac 1 {q_j+q_{j+1}} \le d(q_j\la,\Z) =
      \abs{q_j\la-p_j} < \frac 1 {q_{j+1}}.
    \]
    We thus get
    $\frac 1 {2q_{j+1}} \le \Om(q_j) \le \frac 1 {q_{j+1}}$, which we plug in~\eqref{eqaltcBq}:
    \[
      \sum_{j=1}^{\infty} \left(\frac{1}{q_j}-\frac{1}{q_{j+1}}\right) \log q_{j+1}
      \le \cB(q) \le
      \sum_{j=1}^{\infty} \left(\frac{1}{q_j}-\frac{1}{q_{j+1}}\right) \log(2 q_{j+1}).
    \]
    The \rhs\ and the \lhs\ differ by exactly
    $\frac1{q_0}\log2=\log2$, and the conclusion follows from the fact
    that, since $q_{j+1}\ge\big(\frac{\sqrt5+1}{2})^j$ for all~$j$,
    the series $\sum \frac{1}{q_{j+1}}\log q_{j+1}$ is bounded by a
    universal constant.    
  \end{proof}

  In dimension~$1$, the Bruno function~\eqref{eqdefBrfcn} is commonly
  utilized, rather than $\cB(q)$, but this makes no difference for the
  results about the radius of convergence of the linearization.
Yoccoz was the first to find a lower bound of the order of $\ee^{-\cB(q)}$
instead of $\ee^{-2\cB(q)}$, as part of his celebrated work \cite{Yoccoz}
(which also establishes the optimality of Bruno's condition for
one-dimensional holomorphic germs, but the methods of the present
paper do not seem suitable to recover such a result).
Our article \cite{FMS18}, which was focused on one-dimensional germs,
derived a similar lower bound from Ecalle's method, but relying on an
arithmetical lemma due to Davie, as in \cite{CarlettiM}, whereas in
this paper we will directly show how the Bruno series appears from the
combinatorics.



\section{Formal linearization and Ecalle's tree formalism}\label{sec:formal}

The first key to Ecalle's tree formalism is to adopt an operator point
of view for diffeomorphisms and vector fields, \ie to consider these
geometric objects (or their formal counterparts) respectively as
algebra automorphisms or derivations acting on the algebra of formal
series in~$\nu$ variables.

\black{Note that, classically, in differential geometry,
the notion of Lie derivative allows one 
to identify the vector fields
of a given smooth manifold~$M$ with the $\R$-linear derivations acting
on the algebra $C^\infty(M)$.
As for a transformation $T\col M\to M$, Koopman's
idea \cite{Koop} of associating with it the operator
$\ph\mapsto \ph\circ T$ has proved very fruitful in ergodic
theory}.\footnote{\black{In
the context of measurable dynamical systems, this is the notion of
operator adjoint to the transformation~$T$ (acting on the space of
measurable complex functions) --- see \eg \cite{CFS}.}}
We will be dealing with the formal counterpart of this.

\subsection{The operator point of view}

We recall that the \emph{total order} of a non-zero formal series $\ph(z) =
\sum_{m\in\N^\nu} \ph_m z^m$ is defined as
the minimal value of~$\abs{m}$ over the multi-integers~$m$ such that
$\ph_m \neq 0$. We denote it by $\ord\ph$.
By convention $\ord 0 = \infty$.

Notice that the formal series without constant term are exactly those
for which the total order is $\ge1$. 
For $k\geq 0$
we denote by~$\Cznu_{\geq k}$ the set of formal power series of order
greater than or equal to $k$. Note that $\tgPnu=\big(\Cznu_{\geq 2}\big)^{\nu}$.

We call \emph{operators} the $\C$-linear endomorphisms of \black{the vector space} $\Cznu$.
An operator~$\Th$ is said to be \emph{tangent-to-identity} if
$\Th-\ID$ increases the total order by at least one unit, \ie
$\ord(\Th\ph-\ph) \ge \ord\ph + 1$ for all $\ph \in \Cznu$.
With any $\nu$-tuple $v=(v_1,\ldots,v_\nu)\in (\Cznu_{\geq 1})^{\nu}$,
\black{\ie any substituable tuple of formal series,} we associate the
\emph{composition operator}
\begla
C_v \col \ph \mapsto \ph\circ v
\edla
\black{(formal analog of the Koopman operator).}
%

\begin{lem}   \label{lemAlgAutom}
For any $\nu\in\N^*$, the map $v\mapsto C_v$ induces a bijection
between~$\tgGnuid$ and the set of all tangent-to-identity algebra
automorphisms of $\Cznu$.
\end{lem}

\begin{proof}
  Let $v=(v_1,\dots,v_\nu) \in \tgGnuid$ and
  $a=v-\tmop{id}\in \tgPnu$. For any formal series $\varphi$, using
  the formal Taylor expansion at $z$ of
  $C_v \varphi(z)=\varphi(v(z))=\varphi(z+a(z))$, one easily checks
  that the operator~$C_v$ is a tangent-to-identity algebra
  automorphism. \black{Since $v_i= \Th z_i$ for each~$i$}, the map
  $v\mapsto C_v$ is injective.

  Let $\Th$ be a tangent-to-identity algebra endomorphism of
  $\Cznu$. For $i\in [\nu]$, we let it act on the formal series~$z_i$
  and define $v_i\defeq \Th z_i$;
  as $\ord(\Th z_i-z_i) \ge \ord z_i + 1 =2$, we have
  $v=(v_1,\dots,v_\nu) \in \tgGnuid$ and it only
  remains to prove that $\Th=C_v$. 

  As $\Th$ and $C_v$ are algebra endomorphisms whose images coincide
  on the formal series $z_1,\dots,z_\nu$, their actions coincide on
  polynomials of $z_1,\dots,z_\nu$ too.
  Let $\ph \in \Cznu$:  for any
  $N\in \N$, we can write $\ph=P_N+\ph_N$ with $P_N$ polynomial
  and $\ph_N\in \Cznu_{\geq N}$, thus
\[
  \ord((\Th -C_v)\ph)= \ord((\Th-C_v)\ph_N)
  \geq \min\big\{ \ord(\Th\ph_N-\ph_N), \ord(C_v\ph_N-\ph_N) \big\}
  > N.
\]
Since this holds for all integers $N$, one can deduce that $(\Th -C_v)\ph=0$ and thus $\Th=C_v$.
\end{proof}


\begin{rem}
  \black{One can show that, more generally, $v\mapsto C_v$ is a
  bijection between $(\Cznu_{\geq 1})^{\nu}$ and the set of all
  algebra endomorphisms of $\Cznu$ (obviously a monoid
  antihomomorphism).}
\end{rem}


Any $\nu$-tuple $v=(v_1,\ldots,v_\nu)\in (\Cznu)^{\nu}$ can
also be considered as a formal vector field, for which we denote the
associated $\C$-linear derivation of $\Cznu$ by
\begla
X_v : \varphi \mapsto \sum_{i=1}^\nu v_i \frac{\partial
  \varphi}{\partial z_i}.
\edla
One easily gets an analog of the previous lemma:
\begin{lem}   \label{lemDeriv}
For any $\nu\in\N^*$, the map $v\mapsto X_v$ induces a bijection
between $(\Cznu_{\geq 1})^{\nu}$ and the set of all $\C$-linear derivations of
$\Cznu$ that do not decrease the total order.
\end{lem}

\begin{proof}
The proof is similar to the previous one.
\end{proof}

We will solve the formal linearization problems of
Sections~\ref{seclinpbgr} and~\ref{seclinpbvf} by finding suitable
expressions for the composition operators~$C_h$ associated with the
sought linearizing formal diffeomorphisms~$h$.

\begin{rem}
  One can easily check that, given $a\in \tgPnu$, the linearization
  equation~(\ref{eqqconj}) for a non-resonant formal
  diffeomorphism $g = R_q \circ (\tmop{id}+a)$ 
  is equivalent to
\begin{equation}\label{eq:linqoperator}
C_{R_q} C_h C_{R_q}^{-1} =  C_h C_{\tmop{id}+a},
\end{equation}
and the linearization equation~(\ref{eq:vecConj}) for a non-resonant formal vector field
$V=V_{\lambda}+a$ is equivalent to
\begin{equation}\label{eq:linloperator}
  [X_{V_\la} , C_h] 
  =C_h X_a,
\end{equation}
using $[\,\cdot,\cdot\,]$ to denote operator commutator.
\end{rem}

In both cases, 
$C_h$ will be obtained as a formally convergent series of
\emph{homogeneous} operators.

\begin{Def}
Let $n\in\Z^\nu$. An operator~$D$ is said to be
\emph{$n$-homogeneous} if, for every $m\in\N^\nu$, $D(\C z^m) \subset \C z^{m+n}$.
\end{Def}

For instance, for each $i\in [\nu]$, the partial derivative operator
$\pa_i \defeq \dfrac{\pa\;}{\pa z_i}$ is $n$-homogeneous with $n=-e_i$, where
\[
e_i \defeq \underset{\substack{\uparrow \\ \text{$i$th position}}}{(0,\ldots,0,1,0,\ldots,0)},
\]
while the operators $z_1\pa_1, \ldots,z_\nu\pa_\nu$ are
$0$-homogeneous, \black{and so is}
\[
  \black{X_{V_\la} = \sum_{i\in[\nu]} \la_i z_i \pa_i.}
  \]
Homogeneity will help us to manipulate formally convergent series of
operators.
First, the total order of formal power series allows us to define
summable families for the formal convergence: a family
$(\phi_i)_{i\in I}$ of formal series is summable if, for any
$k\in\N$, the set $\Set{ i\in I \given \ord \phi_i \leq k }$ is finite.
Correspondingly, we can define pointwise summability of families of 
$\C$-linear endomorphisms of $\Cznu$.
The following elementary summability criterion will be of constant use:

\begin{lem}\label{lem:somfam}
  Let $(D_j)_{j\in J}$ be a family of operators such that
  \begin{enumerate}[(i)]
  \item for each $j \in J$, there exists $n_j\in\Z^\nu$ such that
    $\abs{n_j}\ge 0$ and $D_j$ is $n_j$-homogeneous;
\item for each $k\in \N$, the set $\Set{j \in J \given \abs{n_j}=k }$
  is finite.
\end{enumerate}
Then, for any $\varphi\in \Cznu$, the family of formal series
$(D_j \varphi)_{j\in J}$ is summable and 
$\ph\mapsto\sum_{j\in J} D_j\ph$ is a well-defined operator of $\Cznu$.
\end{lem}

\subsection{\black{The set $\cN \subset \Z^\nu$ of homogeneity degrees of the
  nonlinear part}}

We now introduce the set of homogeneity degrees with which we will
deal in practice.
%
Let 
\beglab{eqdefcN}
\cN \defeq \Set[\big]{ m-e_i \given m \in \N^\nu, \;\; \abs{m} \ge 2, \;\; i \in [\nu] }
\subset \Z^\nu,
\edla
so that, either dealing with $f=\tmop{id}+a \in\tgGnuid$ or a vector
field $V=V_{\lambda}+a$, the nonlinear part $a\in \tgPnu$ can be
written
\beglab{eqfwithcN}
a = (a_1,\ldots,a_\nu), \qquad
a_i(z_1,\ldots,z_\nu) = \sum_{n\in\cN} a_{i,n} z^{n+e_i}
\edla
%
with an arbitrary collection of complex numbers
$(a_{i,n})_{(i,n)\in[\nu]\times\cN}$ such that\footnote{\black{These
    are the same coefficients as in~\eqref{eqgenftgGnuid} or~\eqref{eqcoefftiaiminVi}, but with a new
    indexation, which will turn out to be advantageous.}}
\beglab{eqfwithcNa}
a_{i,n} \neq0 \quad\Longrightarrow\quad n+e_i\in\N^\nu.
\edla
%

Observe that $z^{n+e_i} = z_i z^n$, but this does not necessarily mean
that~$a_i$ is divisible by~$z_i$ in the ring $\Cznu$, because a multi-integer
$n\in\cN$ may have a negative component (but at most one, and then the
value of this component must be $-1$).

  \black{The nonlinear part $a \in \tgPnu$ gives rise to a derivation that can be
    decomposed into homogeneous components as follows:}
    \beglab{eqdecompahomog}
    \black{X_a=\sum_{i\in[\nu]} a_i(z) \pa_i = \sum_{n\in\cN} E_n,} \quad
      \black{E_n \defeq \sum_{i\in[\nu]} a_{i,n} z^{n+e_i} \pa_i
      \ens\;
      \text{$n$-homogeneous for each $n\in\cN$.}}
      \edla


      \color{black}
      
  \begin{rem} \label{remNRcN}
    Let us define, for any $q\in (\C^*)^\nu$ and $\la \in \C^\nu$,
    \begla
    \De_q(n)\defeq q^n- 1, \quad \ti\De_\la(n)\defeq\lambda \cdot n
    \qquad \text{for all $n\in\Z^\nu$.}
\edla
        The non-resonance condition for diffeomorphisms can be rephrased as
\beglab{eqqnonrescN}
\text{$q$ non-resonant} \quad\Longleftrightarrow\quad
\text{$\De_q(n)\neq0$ for all $n\in\cN$}
\edla
whereas, for vector fields,
\beglab{eqlanonrescN}
\text{$\lambda$ non-resonant} \quad\Longleftrightarrow\quad
\text{$\ti\De_\la(n)\neq 0$ for all $n\in\cN$.}
\edla
The conjugacy equation~\eqref{eqqconj} with $f=\id+a$ can be
  written
\beglab{eqqconji}
L_q\ppi h_i = a_i\circ h \quad \text{for $i=1,\ldots,\nu$,}
\edla
with operators $L_q\ppi \col \ph \mapsto q_i\ii \ph\circ R_q - \ph$,
  of which the numbers $\De_q(n)$ are eigenvalues:
\beglab{eqDeqnev}
L_q\ppi z^{n+e_i} = \De_q(n) z^{n+e_i}
\quad \text{for all $n\in\Z^\nu$ and $i\in[\nu]$ such that $n+e_i\in\N^\nu$.}
\edla
The conjugacy equation~\eqref{eq:vecConj} with
  $V_i(z)=\la_iz_i+a_i(z)$ for each~$i$ can be written
\beglab{eqlaconji}
\ti L_\la\ppi h_i = a_i\circ h  \quad \text{for $i=1,\ldots,\nu$,}
\edla
with operators $\ti L_\la\ppi \defeq X_{V_\la} - \la_i \id$,
  of which the numbers $\ti\De_\la(n)$ are eigenvalues:
\beglab{eqDelanev}
\ti L_\la\ppi z^{n+e_i} = \ti\De_\la
(n) z^{n+e_i}
\quad \text{for all $n\in\Z^\nu$ and $i\in[\nu]$ such that $n+e_i\in\N^\nu$.}
\edla
\end{rem}

\color{black}

\subsection{\black{$\cN$-forests and $\cN$-trees}}
\label{seccNFT}

\black{We want to use the set~$\cN$ to define the $\cN$-forests and
$\cN$-trees that will index summable families of homogeneous operators
leading to the solution to the formal linearization problems.}
We thus briefly review the relevant formalism following
\cite{FMS18}, to which we refer for more details.

Recall that an \emph{arborescent poset} is a finite poset each element
of which has at most one direct predecessor;
any $\cN$-forest~$F$ can be viewed as an
\emph{$\cN$-decorated arborescent poset}, \ie an arborescent poset
$(V_F,\, \preceq_F)$ endowed with a decoration map
$N_F\col V_F\to\cN$, up to isomorphism---see
Sec.~3 of~\cite{FMS18}. This definition becomes clear when
drawing (decorated) Hasse diagrams:

Given $(n_1,\dots,n_6)\in \cN^6$, the $\cN$-forest 
\begin{equation} \label{eqexempleNF}
F = \BeunedeuxSetroisequatre \quad \Becinqesix
\end{equation}
can be represented by a six-element poset, with vertices
$V_F=\{a,b,c,d,e,f\}$, arborescent order
\[
a \preceq_F b,\ a\preceq_F c\preceq_F d,\ e\preceq_F f
\]
and decoration map
\[
N_F(a)=n_1,\  N_F(b)=n_2,\ N_F(c)=n_3,\ N_F(d)=n_4,\ N_F(e)=n_5,\ N_F(f)=n_6.
\]
Rather than recalling the definitions in \cite{FMS18}, let us
illustrate them on this example:
\begin{itemize}
\item
  The family of all $\cN$-forests, including the empty $\cN$-forest $\emptyset$,
  is denoted by $\cF(\cN)$.
\item The \emph{size} of an $\cN$-forest~$F$, denoted by $\#F$, is
  defined as the cardinality of the underlying set of vertices~$V_F$
  (in the example $\#F=6$; of course, $\#\emptyset=0$).
\item
  The \emph{degree} of~$F$, denoted by~$\deg(F)\in\N$, is the number
  of minimal elements (sometimes called roots) of $(V_F,\, \preceq_F)$ (in the example
  $\deg(F)=2$, and $\deg(\emptyset)=0$).
\item An $\cN$-tree is an $\cN$-forest of degree~$1$; the set of
  $\cN$-trees is denoted by $\cT(\cN)$ (beware that the empty $\cN$-forest
  $\emptyset$ is not an $\cN$-tree).
\item
The \emph{height} $h(F)$ of~$F$ is the maximal cardinality of a
chain\footnote{%
A chain of a poset is a subposet which is totally ordered.
}
of $(V_F,\, \preceq_F)$ (in the example $h(F)=3$ : $a\preceq_F c
\preceq_F d$).
\item If $\sigma\in V_F$, then we denote by $\Tree(\sig,F)$ the $\cN$-subtree
of~$F$ rooted at~$\sig$
(whose set of vertices is $\Set{ \mu \in V_F \given \sig \preceq_F \mu
}$, with the arborescent poset structure induced by~$\preceq_F$).
In our example, $\Tree(c,F)=\Bedeuxtrois$.
\end{itemize}

The set $\cF(\cN)$ can be considered as the free commutative monoid on
$\cT(\cN)$ (with the empty $\cN$-forest $\emptyset$ as the unit). In our
example, we can write~$F$ as the product
$F=T_1 T_2=T_2 T_1$ with
\[ T_1= \BeunedeuxSetroisequatre \quad\text{and}\quad T_2=\Becinqesix
\]
and any $F\in\cF(\cN)$ can be written
\begin{equation}\label{eqdecompF}
F = T_1^{d_1} \cdots T_s^{d_s}, \qquad 
T_1,\ldots,T_s \in \cT(\cN) \ens\text{pairwise distinct}, \quad
d_1,\ldots,d_s \in \N^*,
\end{equation}
a decomposition which is unique up to a permutation of the pairs
$(T_i,d_i)$.
The sum of the
multiplicities is the degree of the $\cN$-forest: $\deg(F) =
d_1+\cdots+d_s$.

\black{Any $\cN$-tree can be written (in a unique way) $T=n \att F$
  with $n\in\cN$ and $F\in\cF(\cN)$, a notation that indicates
  that~$T$ is obtained by attaching a root bearing the decoration~$n$ to
the $\cN$-forest~$F$.} For instance:
\[
  %
%
  n\att 
    \ppb{\BeunedeuxSetroisequatre \ens \Becinqesix}
  =
  \BeBeunedeuxSetroisequatreSecinqesix
\]

$\cN$-trees are $\cN$-forests of degree~$1$.
``Bamboo trees'' correspond to the case when the poset is totally
ordered; we then use the notation
\beglab{eqdefbamboon}
[n_1,\ldots,n_r] \defeq n_1 \att ( n_2 \att ( \cdots ( n_r \att
\emptyset ) \cdots))
\quad\text{for arbitrary}\ens n_1,\ldots,n_r\in\cN
\edla
\black{(this is the same thing as a word $n_1\cdots n_r$ on the alphabet~$\cN$).}
In particular, $[n_1]$ 
is the one-vertex poset with decoration~$n_1$.

Since~$\cN$ is a subset of the commutative group~$\Z^\nu$, we can
define the weight of an $\cN$-forest~$F$ and the weights of its vertices by adding
decorations:
\[
\norm{F} \defeq \sum_{\sig \in V_F} N_F(\sig) \in \Z^\nu, \qquad
\htx\sig \defeq \norm{\Tree(\sig,F)} \in \Z^\nu \ens\text{for $\sig \in V_F$.}
\]
Finally, we define $\abs{\,\cdot\,}: \cN \rightarrow \N^{*}$ by $|(n_1,\dots ,
n_{\nu})|= n_1 +\dots + n_{\nu}$ and
\beglab{ineqabsnormF}
\abs{\norm{F}} \defeq \sum_{\sig \in V_F} \abs{N_F(\sig)}\geq \#F.
\edla

\subsection{\black{The coarmould associated with the nonlinear part}}

We call coarmould any family of operators indexed by $\cF(\cN)$.
The summation in~\eqref{eqmldexpaCh} will be over all $\cN$-forests
$F\in\cF(\cN)$ and we now define the coarmould $D_\bul=
(D_F)_{F\in\cF(\cN)}$
associated with the nonlinear part $a\in \tgPnu$ of a formal
diffeomorphism $f=\tmop{id}+a \in\tgGnuid$ or vector field
$V=V_{\lambda}+a$.
Recall that we use~$\pa_i$ as shorthand for $\dfrac{\pa\;}{\pa z_i}$.
\begin{lem}    \label{lemdefDbulf}
Let $a \in \tgPnu$ be defined by coefficients $(a_{i,n})_{(i,n)\in[\nu]\times\cN}$ as in~\eqref{eqfwithcN}--\eqref{eqfwithcNa}.
The formulas \begin{enumerate}
\item[\textup{(i)}]
$D_\emptyset = \ID$,
\item[\textup{(ii)}]
$D_{n\att F} = \big(D_F(a_{1,n} z^{n+e_1})\big)\pa_1 + \cdots +
\big(D_F(a_{\nu,n} z^{n+e_\nu})\big)\pa_\nu\;$
for any $n\in\cN$ and $F \in \cF(\cN)$,
\item[\textup{(iii)}]
$\dst D_{U_1\cdots U_d} = \frac{1}{d_1!\cdots d_r!} \sum_{i \in [\nu]^d}
(D_{U_1} z_{i_1}) \cdots (D_{U_d} z_{i_d})
\pa_{i_1} \cdots \pa_{i_d}\;$
for any $U_1,\ldots, U_d \in \cT(\cN)$, with $d\ge1$,
where $r$ is the number of distinct $\cN$-trees in this list and
$d_1,\ldots,d_r$ are the corresponding multiplicities
(\ie $U_1\cdots U_d$ can be written $T_1^{d_1}\cdots T_r^{d_r}$ as in~\eqref{eqdecompF})
\end{enumerate}
determine uniquely a coarmould $D_\bul$.
Moreover, $D_F$ is $\norm{F}$-homogeneous for each $F \in\cF(\cN)$.
\end{lem}


\begin{proof}
  We first notice that, although the decomposition in disjoint union of
  $\cN$-trees of an $\cN$-forest $F = U_1\cdots U_d$ is determined
  only up to a permutation of $(U_1,\ldots, U_d)$, the formula in~(iii)
  makes sense because the \rhs\ is invariant under permutation.
  We also observe that~(iii) is compatible with~(ii) in the sense
  that, when $U=n\att F \in\cT(\cN)$,~(iii) merely asserts that
  $D_U = (D_U z_1)\pa_1 + \cdots + (D_U z_\nu)\pa_\nu$,
  \ie that~$D_U$ is a derivation, while~(ii) tells which
  derivation it is in terms of~$n$ and $D_{F}$.
 The operator~$D_F$ is  well-defined for any $F\in \cF(\cN)$
 by induction on the height of~$F$ and its homogeneity property is also easily obtained by the same induction.
\end{proof}


\begin{Def}
  For any $a \in \tgPnu$, the coarmould defined in
  Lemma~\ref{lemdefDbulf} is called the \emph{coarmould associated
    with~$a$} and denoted by~$D_\bul(a)$, or simply~$D_\bul$ if~$a$ is
  clear from the context.
\end{Def}

\black{Notice that, in the case of size~$1$,
  Lemma~\ref{lemdefDbulf}(ii) yields
  $D_{[n_1]}=E_{n_1}$, the $n_1$-homogeneous component of the
  derivation associated with~$a$ as in~\eqref{eqdecompahomog}.}

The above definition of~$D_\bul$ is the natural
$\nu$-dimensional extension of the definition given in \cite{FMS18} for the one-dimensional case.
The concrete shape of~$D_F$ given in \cite{FMS18} can also be extended to dimension~$\nu$:
\begin{lem}   \label{lemDFbetaiF}
For any finite set~$R$ and any map $i \col R \to [\nu]$,
we define a $0$-homogeneous operator
\[
\de_i \defeq \bigg( \prod_{\rho\in R} z_{i(\rho)} \bigg) \prod_{\rho\in R} \pa_{i(\rho)}
= z_1^{\di1} \cdots z_\nu^{\di\nu} \pa_1^{\di1} \cdots \pa_\nu^{\di\nu},
\]
where $\di1 \defeq \card\big( i\ii(1) \big)$, \ldots, 
$\di\nu \defeq \card\big( i\ii(\nu) \big)$.
Then, for any $\cN$-forest~$F$, 
choosing a representative $(V_F,\, \preceq_F,\, N_F)$ and denoting
by~$R_F$ the set of all minimal elements of $(V_F,\, \preceq_F)$,
we have
\beglab{eqDFbeiF}
D_F(a) = z^{\norm{F}} \!\sum_{i \colb R_F \to [\nu]} \be_{i,F}(a) \, \de_i,
\edla
where the scalar coefficients $\be_{i,F}(a)$ are polynomial functions
with rational non-negative coefficients in the variables
$(a_{i,n})_{(i,n)\in[\nu]\times\cN}$.
\end{lem}
\begin{proof}
Induction on the height of~$F$.
\end{proof}


For instance, for $\nu = 2$ and $\deg(F) = d \le 2$, choosing the
representative $(V_F,\, \preceq_F,\, N_F)$ so that $R_F =
\{1,\ldots,d\}$ and identifying any map $i\col R_F \to [\nu]$ with the
$d$-tuple $\big( i(1), \ldots, i(d) \big)$,
we get
\begin{align*}
\deg(F) &= 1 &\Longrightarrow&&
D_F &= z^{\norm{F}} 
\big( \be_{(1),F} z_1 \pa_1 + \be_{(2),F} z_2 \pa_2 \big) \\[1ex]
\deg(F) &= 2 &\Longrightarrow&&
D_F &= z^{\norm{F}} 
\big( \be_{(1,1),F} z_1^2 \pa_1^2 + (\be_{(1,2),F} + \be_{(2,1),F}) z_1z_2 \pa_1\pa_2
+ \be_{(2,2),F} z_2^2 \pa_2^2 \big).
\end{align*}
The precise value of the coefficients~$\be_{i,F}$ is given in
Appendix~\ref{secappcomputDF} in full generality, but this information
will not be used in the paper.
What we \emph{will} use is 
(to implement a majorant series argument)
the fact that the coefficients of the polynomials $\be_{i,F}(a)$ are
non-negative---see~\eqref{ineqmajorantDT}.


The key algebraic property of the operators $D_F=D_F(a)$ is what
Ecalle calls coseparativity:
%
the way they act on products of formal power series is related to the
product law of the monoid~$\cF(\cN)$.

\begin{prop}   \label{propcosepcontractsepnu}
The coarmould~$D_\bul$ associated with any given $a \in \tgPnu$ satisfies
\beglab{eqcoloiDF}
\dst D_F(\ph\psi) = 
\sum_{\substack{ (F',F'') \in \cF(\cN)\times\cF(\cN) \\
           \text{such that}\ F = F' F''}}
(D_{F'}\ph) (D_{F''}\psi)
\edla
for any $F\in\cF(\cN)$ and $\ph,\psi \in \Cznu$.
\end{prop}


\begin{proof}
Let us prove this statement for arbitrary $F\in\cF(\cN)$ and
$\ph,\psi \in \Cznu$.
If $F=\est$ then \eqref{eqcoloiDF} is obvious.
We thus suppose
$F = U_1 \cdots U_d = T_1^{d_1} \cdots T_r^{d_r}$ as in point~(iii) of
Lemma~\ref{lemdefDbulf}.
More precisely, we assume that $d\ge1$ and
we are given a partition $[d] = P_1 \sqcup \cdots \sqcup P_r$
and collections of
$\cN$-trees $(U_1, \ldots, U_d)$ and $(T_1,\ldots,T_r)$ such that, for each $k\in[r]$,
\[
d_k = \card P_k \ge 1, \qquad\qquad
j\in P_k \imp T_k = U_j.
\]
We have
\[
D_F = \frac{1}{d_1!\cdots d_r!} \sum_{i \in [\nu]^d}
(D_{U_1} z_{i_1}) \cdots (D_{U_d} z_{i_d})
\pa_{i_1} \cdots \pa_{i_d}
\]
and, for each $i \in [\nu]^d$, the Leibniz rule yields
\[
\pa_{i_1} \cdots \pa_{i_d}(\ph\psi) =
\sum_{\De\subset[d]} \bigg(\prod_{j\in \De}\pa_{i_j}\bigg)\ph
\cdot \bigg(\prod_{j\in \dde}\pa_{i_j}\bigg)\psi,
\]
with summation over all subsets of~$[d]$ including
the empty set, hence
\begin{multline*}
D_F(\ph\psi) = \frac{1}{d_1!\cdots d_r!} \sum_{\De\subset[d]} \chi_\De 
\quad\text{with}\\[1ex]
\chi_\De \defeq
\sum_{i \in [\nu]^d}
\bigg(\prod_{j\in \De} D_{U_j} z_{i_j} \bigg)
\bigg(\prod_{j\in \De}\pa_{i_j}\bigg)\ph
\cdot \bigg(\prod_{j\in \dde} D_{U_j} z_{i_j} \bigg)
\bigg(\prod_{j\in \dde}\pa_{i_j}\bigg)\psi.
\end{multline*}
Now, for each $\De\subset[d]$, there is a
bijection $[\nu]^d \to [\nu]^{\De} \times [\nu]^{\dde}$ obtained by
viewing each $i \in [\nu]^d$ as a map $[d]\to[\nu]$ and sending it to
the pair of restrictions $(i',i'')$, $i' \defeq i_{|\De}$, $i'' \defeq
i_{|\dde}$, thus $\chi_\De = \ph_{\De} \cdot\psi_{\dde}$ with
\[
\ph_{\De} =
\sum_{i'\in[\nu]^{\De}} \bigg(\prod_{j\in \De} D_{U_j} z_{i'_j} \bigg)
\bigg(\prod_{j\in \De}\pa_{i'_j}\bigg)\ph, \quad
\psi_{\dde} = \sum_{i''\in[\nu]^{\dde}} \bigg(\prod_{j\in \dde} D_{U_j} z_{i''_j} \bigg)
\bigg(\prod_{j\in \dde}\pa_{i''_j}\bigg)\psi.
\]

Let us map each $\De\subset[d]$ to the pair of $\cN$-forests
$(F',F'')$ such that $F=F'F''$ defined as
\[
F' \defeq \prod_{j\in \De} U_j = T_1^{d'_1} \cdots T_r^{d'_r}, \quad
F'' \defeq \prod_{j\in \dde} U_j = T_1^{d''_1} \cdots T_r^{d''_r},
\]
with $d'_k = \card(\De\cap P_k)$, $d''_k = d_k-d'_k$ for each $k\in[r]$.
This map is surjective and the number of preimages of any such
$(F',F'')$ is $\binom{d_1}{d'_1} \cdots \binom{d_r}{d'_r}$
(indeed, a preimage~$\De$ is determined choosing 
$\De_1 \subset P_1$, \ldots, $\De_r \subset P_r$ and taking their union,
with the sole constraint that $\card \De_1 = d'_1$, \ldots, $\card \De_r = d'_r$).
Since
\[
\ph_{\De} = d'_1! \cdots d'_r! \, D_{F'}\ph, \quad \psi_{\dde} = d''_1! \cdots d''_r! \, D_{F''}\psi,
\]
we get 
\[
D_F(\ph\psi) = \frac{1}{d_1!\cdots d_r!} 
\sum_{\substack{ (F',F'') \\
           \text{s.t.}\, F = F' F''}}
\binom{d_1}{d'_1} \cdots \binom{d_r}{d'_r}
d'_1! \cdots d'_r! \, D_{F'}\ph \cdot d''_1! \cdots d''_r! \, D_{F''}\psi,
\]
which yields the desired conclusion.
\end{proof}

\subsection{Armoulds and tree expansions}

We call armould any family of complex numbers indexed by~$\cF(\cN)$
and use the notation~$A^\bul$ for an armould $(A^F)_{F\in\cF(\cN)}$.
One easily checks that the coarmould associated to a given $a \in
\tgPnu$ gives rise to a pointwise summable family of operators
$(A^F D_F(a))_{F\in\cF(\cN)}$ for any armould~$A^\bul$.
Indeed, the summability criterion of Lemma~\ref{lem:somfam} is
fulfilled because each $D_F(a)$ is $\norm{F}$-homogeneous with $\abs{\norm{F}}\ge0$
and, for every $k\ge0$, there are only finitely many
$\cN$-forests~$F$ such that $\abs{\norm{F}}\leq k$
(because of~\eqref{ineqabsnormF} and~\eqref{eqdefcN}).
%
%
We thus have 
\begin{lem}
  Let $a\in\tgPnu$.
%
%
For any armould~$A^\bul$, the series
\beglab{eqdeftreeexp}
\sum A^\bul D_\bul(a) \defeq \sum_{F\in\cF(\cN)} A^F D_F(a)
\edla
defines an operator of $\Cznu$.
\end{lem}

A series of the form~\eqref{eqdeftreeexp} is called a tree expansion.
Dually to Proposition~\ref{propcosepcontractsepnu}, we may give
simple conditions on the armould~$A^\bul$ to ensure that the tree
expansion~\eqref{eqdeftreeexp} is a derivation or a
tangent-to-identity algebra automorphism of $\Cznu$:

\begin{prop} \label{prop:arex} Let $a \in \tgPnu$ and $D_\bul$ its associated coarmould. 
\begin{itemize}
\item[\textup{(i)}] If an armould $A^\bul$ is separative, that is
  $A^\emptyset=1$ and, for any $d\ge1$,
  \begla
  U_1,\ldots, U_d \in \cT(\cN) {\quad\Longrightarrow\quad}
  A^{U_1\cdots U_d}=A^{U_1}\cdots A^{U_d},
  \edla
  then the tree expansion $\sum A^\bul D_\bul$ coincides with the
  composition operator~$C_\th$ where $\th = (\th_1,\ldots,\th_\nu)
  \in \tgGnuid$ is defined by
\beglab{eqdefthiz}
\th_i(z) = z_i + \sum_{T \in \cT(\cN)} A^T D_T z_i
= z_i + \sum_{T \in \cT(\cN)} A^T \be_{i,T} \, z^{\norm T + e_i}
\quad \text{for} \ens i=1,\ldots,\nu.
\edla
%
%
\item[\textup{(ii)}] If an armould $A^\bul$ is antiseparative,
  that is $A^\emptyset=0$ and 
  \begla
  U_1,\ldots, U_d \in \cT(\cN) \ens\text{with}\ens d\ge2 {\quad\Longrightarrow\quad}
  A^{U_1\cdots U_d}=0,
  \edla
  then the tree expansion $\sum A^\bul D_\bul$ coincides with the
  derivation $X_v$ where $v = (v_1,\ldots,v_\nu)\in (\Cznu_{\geq
    1})^{\nu}$ is defined by $\dst v_i\defeq \sum_{T \in \cT(\cN)} A^T D_T
  z_i$ for $i=1,\ldots,\nu$.
\end{itemize}

\end{prop}

\begin{proof}
(i) Suppose that $A^\bul$ is a separative armould.
Since $A^\est D_\est=\ID$ and $D_F(a)$ is $\norm{F}$-homogeneous with
$\abs{\norm{F}}\ge1$ for $F\neq\est$,
we see that $C \defeq \sum A^\bul D_\bul$ is a
tangent-to-identity operator.
Given $\ph,\psi \in \Cznu$, thanks to Proposition~\ref{propcosepcontractsepnu},
\[
  C(\ph \psi) = \sum_{F \in  \cF(\cN)} A^F D_F(\ph \psi) =
  \sum_{F \in  \cF(\cN)} \, \sum_{\substack{ (F',F'') \in \cF(\cN)\times\cF(\cN) \\
           \text{such that}\ F = F' F''}}
       A^{F'F''}(D_{F'}\ph) (D_{F''}\psi).
     \]
Separativity implies that $A^{F'F''}=A^{F'}A^{F''}$, thus 
\begin{align*}
  C(\ph \psi)&=\sum_{F \in  \cF(\cN)} \, \sum_{\substack{ (F',F'') \in \cF(\cN)\times\cF(\cN) \\
           \text{such that}\ F = F' F''}}
A^{F'}A^{F''}(D_{F'}\ph) (D_{F''}\psi) \\[.8ex]
             &=\sum_{(F',F'') \in \cF(\cN)\times\cF(\cN)}(A^{F'}D_{F'}\ph) (A^{F''}D_{F''}\psi)
= C(\ph)C(\psi).
\end{align*}
This shows that~$C$ is an algebra automorphism and
%
Lemma~\ref{lemAlgAutom} thus yields $C=C_\th$ with
$\th = (Cz_1,\ldots,Cz_\nu)$.
Formula~\eqref{eqdefthiz} follows from
\begla
\deg(F) \ge 2 \quad\Longrightarrow\quad D_F z_i=0 \quad\text{for}\ens i=1,\ldots,\nu
\edla
(consequence of Lemma~\ref{lemdefDbulf}(iii)).
\medskip

\noindent (ii)
In view of Lemma~\ref{lemDeriv}, the proof of the second point is even simpler since the tree expansion is a series of operators indexed by trees that are, by definition, derivations.
\end{proof}



Here are the simplest possible examples:
\begin{itemize}
\item
  The armould~$I^\bul$ defined by
  $I^F\defeq 1$ if the height of~$F$ is at most~$1$ and~$0$ otherwise
  %
  is separative, and
  \begla
\sum I^\bul D_\bul(a)=C_{\tmop{id}+a}.
\edla
\item
  The armould~$J^\bul$ defined by
  $J^F\defeq 1$ if $\#F=1$ and~$0$ otherwise
  %
  is antiseparative, and
  \begla
\sum J^\bul D_\bul(a) = X_a. 
\edla
\end{itemize}

These formulas are obvious consequences of Proposition~\ref{prop:arex},
using the identity
\begla
D_{[n]} z_i = a_{i,n} z^{n+e_i}
\quad \text{for all} \ens
n\in\cN
\ens\text{and}\ens
i\in[\nu]
\edla
(with notation~\eqref{eqdefbamboon}), which follows from Lemma~\ref{lemdefDbulf}(ii).

Let us conclude this section with a result on``geometric''
separative armoulds that will be useful when dealing with the
analyticity of diffeomophisms associated to tree expansions.

\begin{prop} \label{prop:geomexpansion}
  Let $a\in \tgPnu$ and $A,B\in\C$ with $B\not=0$.  
  Then the
  armould $K^\bul$ defined by:
  \begin{equation}
    %
    %
    K^F \defeq A^{\#F}B^{\abs{\norm{F}}}
    \quad \text{for}\ens F \in \cF(\cN)
\end{equation}
is separative and the corresponding tree expansion is
\begla
\sum K^\bul D_\bul(a) = C_w
\edla
where $w = (w_1,\ldots,w_\nu)$ is the composition inverse in
$\tgGnuid$ of the tangent-to-identity formal diffeomorphism 
\beglab{eqdefg}
\black{ G(z) \defeq z - A B\ii a(Bz).}
%
%
\edla
\end{prop}

Note that in the case $(A,B)=(-1,1)$ the formal diffeomorphism~$w$ is
nothing but the composition inverse of $\tmop{id}+a$.

\begin{proof}
  The armould $K^\bul$ is clearly separative;
  by Proposition~\ref{prop:arex}(i), its tree expansion coincides with
  a composition operator $C_w$, with a formal tangent-to-identity
  diffeomorphism~$w$ whose components are
  $w_i \defeq z_i +\sum_{T \in \cT(\cN)} K^T D_T z_i$.
  We must check that this~$w$ is indeed the composition inverse of~$G$
  of~\eqref{eqdefg}.
  %

  \color{black} We note that~$K^\bul$ solves the elementary ``armould equation''
  \[
    K^{n\att F} = A B^{\abs n} K^F
    \quad \text{for all $F\in \cF(\cN)$ and $n\in\cN$.}
  \]
Since any $\cN$-tree can be uniquely written as $T=n\att F$ with
$F \in \cF(\cN)$ and $n\in\cN$, formula~\eqref{eqdefthiz} yields
  \begin{align*}
    w_i(z) = z_i + \sum_{T\in\cT(\cN)} K^T D_T z_i
    &= z_i + \sum_{F \in \cF(\cN)} \, \sum_{n\in\cN} K^{n\att F}D_{n\att F} z_i \\[.8ex]
 &= z_i + \sum_{F \in \cF(\cN)} \, \sum_{n\in\cN} A B^{\abs n} K^F D_F(a_{i,n} z^{n+e_i}) \\[.8ex]
 &= z_i+ \sum_{F \in \cF(\cN)} K^F D_F\left(\sum_{n\in\cN} A B^{\abs n}
   a_{i,n} z^{ n+e_i}\right),
  \end{align*}
  where the last-but-one step was consequence of
  Lemma~\ref{lemdefDbulf}(ii).
  In the latter expression, the inner sum is nothing but $A B\ii a_i(Bz)$
  and, since $\sum K^F D_F = C_w$, we get
  %
  %
  \[
    w_i(z) = z_i+A B\ii a_i(Bw(z)), 
  \]
  which amounts to $G_i\circ w(z)=z_i$, whence the conclusion
  follows.
  \color{black}
\end{proof}

\subsection{ The set $\cF^+(\cN) \subset \cF(\cN)$ and the armoulds $S^\bul(q)$ and $\tilde{S}^\bul(\lambda)$}

We recall that, for any $F \in \cF(\cN)$,
\[
%
  %
%
  \htx\sig = \norm{\Tree(\sig,F)}
= \sum_{\mu \in V_F \ \text{s.t.}\ \sig \preceq_F \mu} N_F(\mu) 
  \in \Z^\nu \quad\text{for}\ens \sig \in V_F.
\]
\color{black}
As mentioned in Remark~\ref{remNRcN}, the non-resonance conditions for
a spectrum~$q$ or~$\la$
guarantee that the complex numbers
$\De_q(n)=q^n - 1$ or $\ti\De_\la(n)=\lambda \cdot n$ do not vanish
for any $n\in\cN$.
\color{black}
However, they may vanish for some nonzero $n\in\Z^\nu$ (what Ecalle calls
``fictitious resonance'' in \cite{EcNonAbord}).
\begin{Def}   \label{DefSbulqnu}
  Let
  \beglab{eqdefcFp}
  \cF^+(\cN) \defeq \Set{ F \in \cF(\cN) \given \htx\sig \in \cN
  \ \text{for all $\sig\in V_F$} }
  \edla
  with the convention $\est \in \cF^+(\cN)$.
For each non-resonant $q\in(\C^*)^\nu$, we define an armould~$S^\bul(q)$ by
\beglab{eqdefSFq}
S^F(q) \defeq \prod_{\sig \in V_F} \black{\frac{1}{\De_q(\shtx\sig)}} 
\quad\text{if $F \in \cF^+(\cN)$,} \qquad
S^F(q) \defeq 0 \quad\text{if $F \in \cF(\cN)\setminus\cF^+(\cN)$.}
\edla
For each non-resonant spectrum $\lambda\in\C^\nu$, we define an armould~$\tilde{S}^\bul(\lambda)$ by
\beglab{eqdeftiSFla}
\tilde{S}^F(\lambda) \defeq \prod_{\sig \in V_F} \black{\frac{1}{\ti\De_\la(\shtx\sig)}} 
\quad\text{if $F \in \cF^+(\cN)$,} \qquad
\tilde{S}^F(\lambda)\defeq 0 \quad\text{if $F \in \cF(\cN)\setminus\cF^+(\cN)$.}
\edla
\end{Def}

\color{black}

In the example~\eqref{eqexempleNF}, the value of the first armould is
\[
  S^F(q) = \frac{1}{\De_q(n_1+n_2+n_3+n_4)  \De_q(n_2)  \De_q(n_3+n_4)
  \De_q(n_5+n_6) \De_q(n_6)}
\]
and the value of the second armould is obtained by replacing~$\De_q$ with~$\ti\De_\la$.

\color{black}


\begin{lem}
The armould~$S^\bul(q)$ and $\tilde{S}^\bul(\lambda)$ are separative for every non-resonant~$q$ or $\lambda$.
\end{lem}


\begin{proof}
We have $S^\est(q) = \ti S^\est(\la)= 1$ by definition. Let $F',F'' \in \cF(\cN)$. We
observe that
\beglab{eqsubmonoidcompl}
F'F'' \in \cF^+(\cN) {\quad\Longleftrightarrow\quad}
\text{$F'$ and $F'' \in \cF^+(\cN)$}
\edla
(in particular~$\cF^+(\cN)$ is a submonoid of $\cF(\cN)$).
This makes it easy to check the identities
\begla
S^{F'F''}(q) = S^{F'}(q)S^{F''}(q), \qquad
\ti S^{F'F''}(\la) = \ti S^{F'}(\la)\ti S^{F''}(\la)
\edla
for all $F',F''\in \cF(\cN)$.
%
\end{proof}

As we shall see now, these two innocent-looking armoulds give the
solution to the formal linearization problems:
the corresponding tree expansions are the composition operators
associated by the sought formal linearizing transformations~$h$.
%
%
%
\black{Moreover, all the potential accumulation of small denominators is encapsulated in the coefficients
$S^F(q)$ and $\tilde{S}^F(\lambda)$ in an efficients way, which will
make it possible to prove analyticity results under Bruno's condition
in Section~\ref{sec:analytic}.}

\subsection{Solution of the formal linearization problem}

Our approach to the formal linearization problems requires a
``vanishing" lemma whose proof relies on tools developed in Section~\ref{sec:vanish}.


\begin{lem}    \label{propnullDcFplus}
The coarmould associated with any $a\in \tgPnu$ satisfies
\beglab{eqnullDcFplus}
F \in \cF(\cN)\setminus \cF^+(\cN) \Imp D_F=0. 
\edla
\end{lem}


\begin{proof} This is an immediate corollary of
  Proposition~\ref{prop:annul2}(i)---see
  Definition~\ref{defunivvanish} and~\eqref{eqFpincludesNV}.
%
%
%
\end{proof}


We are now ready for the formal part of the main results announced in Section~\ref{intro}.


\begin{thm}    \label{thmformaldiffeonu}
  The composition operator of the solution~$h$ to the formal
  linearization problem for $g = R_q \circ(\tmop{id}+a)$ with
  non-resonant $q \in (\C^*)^\nu$ and $a\in \tgPnu$ is
\begla
C_h = \sum S^\bul(q) D_\bul(a),
\edla
where $S^\bul(q)$ is the separative armould defined in
Definition~\ref{DefSbulqnu} and~$D_\bul(a)$ is the coarmould associated
with~$a$.
The solution itself is $h = (h_1,\ldots,h_\nu) \in \tgGnuid$, with
\beglab{eqdefhievalcontr}
h_i(z) = z_i + \sum_{T\in\cT(\cN)\cap\cF^+(\cN)} \,
\prod_{\sig \in V_T} \frac{1}{q^{\shtx\sig} - 1} 
\, D_T z_i
\quad \text{for}\ens i=1,\ldots,\nu.
\edla
\end{thm}

\begin{thm}    \label{thmformalfieldnu}
  The composition operator of the solution~$h$ to the formal
  linearization problem for a vector field $V =V_\lambda +a$ with
  non-resonant $\lambda \in \C^\nu$ and $a \in \tgPnu$ is
\begla
C_h = \sum \tilde{S}^\bul(\lambda) D_\bul(a),
\edla
where $\tilde{S}^\bul(\lambda)$ is the separative armould defined in
Definition~\ref{DefSbulqnu} and~$D_\bul(a)$ is the coarmould associated
with~$a$.
The solution itself is $h = (h_1,\ldots,h_\nu) \in \tgGnuid$, with
\beglab{eqdefhievalcontrvf}
h_i(z) = z_i + \sum_{T\in\cT(\cN)\cap\cF^+(\cN)} \,
\prod_{\sig \in V_T}  \frac{1}{\lambda \cdot \hat{\sig}} \,
D_T z_i
\quad \text{for}\ens i=1,\ldots,\nu.
\edla
\end{thm}


\color{black}

\begin{rem}
We thus can say that the coefficients~$h_{i,m}$ of~\eqref{eqdefhim}
are given by explicit finite sums:
for each $i\in[\nu]$ and $m\in\N^\nu$ such that $\abs m\ge2$,
\begla
h_{i,m} = \sum_{T} \cS^T \be_{i,T}(a)
\edla
with summation over all $T\in\cT(\cN)\cap\cF^+(\cN)$ such that $\norm
T = m-e_i$,
with $\cS^T=S^T(q)$ or $\tilde{S}^T(\lambda)$
and $\be_{i,T}(a)$ as in Lemma~\ref{lemDFbetaiF} (or, even more explicitly,~\eqref{eqexplicitbeiF}).
\end{rem}

\color{black}


\begin{proof}[Proof of Theorem~\ref{thmformaldiffeonu}]
The definition~\eqref{eqdefSFq} of~$S^\bul(q)$ shows that 
\beglab{eqarmldeqndiffeo}
n\att F \in \cT(\cN)\cap\cF^+(\cN) \Imp
\black{\De_q(n+\norm{F}) S^{n\att F}(q) = S^F(q).}
%
%
\edla
Ultimately, the result follows from the fact that~$S^\bul(q)$ is a
separative armould satisfying~\eqref{eqarmldeqndiffeo}.

Indeed, let $\Th \defeq \sum S^\bul(q) D_\bul$.
Since $S^\bul(q)$ is separative, by
Proposition~\ref{propcosepcontractsepnu} we have $\Th = C_{h^*}$ with
$h^* = (h^*_1, \ldots, h^*_\nu)$ and
\beglab{eqdefhstnu}
h^*_i(z) \defeq \Th z_i = z_i + \sum_{T\in\cT(\cN)} S^T(q) D_T z_i
= z_i + \sum_{T\in\cT(\cN) \cap \cF^+(\cN)} S^T(q) D_T z_i
\edla
(where the last identity results from
Lemma~\ref{propnullDcFplus}), \ie $h^*_i(z)$ coincides with the
\rhs\ of~\eqref{eqdefhievalcontr}.
\color{black}
We must show that this~$h^*$ is a formal solution to the conjugacy
equation, which we write in the form~\eqref{eqqconji} as in Remark~\ref{remNRcN},
%
\ie we must compute the action of $L_q\ppi \col \ph \mapsto q_i\ii \ph\circ R_q - \ph$
on~$h^*_i$ and check that it gives $a_i\circ h^*$ for each $i=1,\ldots,\nu$.

We see that the operator~$L_q\ppi$ annihilates~$z_i$ and, for each
$T\in\cT(\cN)\cap\cF^+(\cN)$, acts on $D_Tz_i$ by multiplication by
$\De_q(\norm T)$, by virtue of~\eqref{eqDeqnev} (note that $D_T
z_i\in\C z^{\norm T+e_i}$ by homogeneity of~$D_T$ and,
since $D_T z_i\in \Cznu$, $\norm T+e_i\in\N^\nu$ or $D_T z_i=0$).
Thus
  \begin{multline*}
    L_q\ppi h^*_i = \sum_{T\in\cT(\cN)\cap\cF^+(\cN)} \De_q(\norm T)
    S^T(q) D_T z_i
    = \sum_{\substack{ (n,F)\in\cN\times\cF(\cN)
 \\ \text{such that}\ n\att F \in\cF^+(\cN) }}
\De_q(n+\norm F) S^{n\att F}(q) D_{n\att F} z_i\\[.5ex]
= \sum_{(n,F)\in\cN\times\cF(\cN)}
S^F(q) D_F(a_{i,n}z^{n+e_i})
= \sum_{F\in\cF(\cN)} S^F(q) D_F  \left( \sum_{n\in\cN} a_{i,n} z^{n+e_i} \right),
    \end{multline*}
    where~\eqref{eqarmldeqndiffeo} was used at the last-but-one step
    together with Lemma~\ref{lemdefDbulf}(ii),
    and we availed ourselves of Lemma~\ref{propnullDcFplus} to drop the requirement $n\att F\in\cF^+(\cN)$.
    The inner sum in the latter expression is~$a_i(z)$, and it is
    acted upon by the operator $\sum S^\bul(q) D_\bul$, which is none
    other than $C_{h^*}$. The result is thus $a_i\circ h^*$, as desired.
\color{black}
\end{proof}

\begin{proof}[Proof of Theorem~\ref{thmformalfieldnu}]
  Ultimately, the result will follow from the fact that~$\ti S^\bul(\lambda)$ is a
  separative armould that satisfies
  \beglab{eqarmldvf}
n \att F \in \cF^+(\cN) \Imp
\black{ \ti\De_\la(n+\norm F)\tilde{S}^{n\att F}(\lambda)
= \tilde{S}^F(\lambda) }
%
%
\edla
(an equation that stems from~\eqref{eqdeftiSFla}).

Indeed, similarly to the previous proof, $\sum \tilde{S}^\bul(\lambda) D_\bul = C_{h^*}$ with
$h^* = (h^*_1, \ldots, h^*_\nu)$ given by
\beglab{eqdefhstnuf}
h^*_i(z) = z_i + \sum_{T\in\cT(\cN)} \tilde{S}^T(\lambda) D_T z_i
= z_i + \sum_{T\in\cT(\cN) \cap \cF^+(\cN)} \tilde{S}^T(\lambda) D_T z_i
\edla
and we just need to check that this~$h^*$ satisfies the conjugacy
equation,
\color{black}
which we write in the form~\eqref{eqlaconji} as in
Remark~\ref{remNRcN}.
With a computation analogous to the one of the previous proof,
replacing the operator $L_q\ppi$ with $\ti L_\la\ppi \defeq X_{V_\la}
- \la_i \id$ and the eigenvalues $\De_q(\norm T)$ with
$\ti\De_\la(\norm T)$ (following~\eqref{eqDelanev}),
we find $\ti L_\la\ppi h^*_i = a_i\circ h^*$.
\color{black}
\end{proof}


\section{Solution of the holomorphic linearization problem}\label{sec:analytic}

In Section~\ref{sec:formal}, the introduction of the set $\cF^+(\cN)$
made it possible to define without ambiguity the ``linearizing
armoulds" $S^\bul(q)$ and $\ti S^\bul(\lambda)$ and obtain the formal
linearization as a sum over this set.
In fact, as will be seen in Section~\ref{secNVcN}, it is possible to
restrict this sum to an even smaller set of $\cN$-forests
$NV(\cN)\subset \cF^+(\cN)$, because 
$F \in \cF^+(\cN) \setminus NV(\cN) \imp D_F(a)=0$ for all $a\in \tgPnu$.
%
If we denote by $\cS^\bullet$ the armould
$S^\bul(q)$ or $\tilde{S}^\bul(\lambda)$,
Theorems~\ref{thmformaldiffeonu} and~\ref{thmformalfieldnu} thus entail
%
that the composition automorphism of the linearizing transformation
for $R_q\circ(\tmop{id}+a)$ or $V_\la+a$ reads
\begin{equation}\label{eqrestrictedseries}
C_h=\sum_{F\in \cF^+(\cN)}\cS^F D_F(a)=\sum_{F\in NV(\cN)}\cS^F D_F(a).
\end{equation}
The definition and main properties of this set of ``non-vanishing"
$\cN$-forests are deferred to Section~\ref{sec:vanish}.

On the other hand, the armould $\cS^\bul=S^\bul(q)$ or $\tilde{S}^\bul(\la)$ depends only
on the spectrum~$q$ or~$\la$, not on the nonlinear part~$a$ of the
dynamical system to be linearized.
We will show in Section~\ref{sec:estim} that,
under the Bruno condition (Definition~\ref{def:Brunodiff} or~\ref{def:BrunoVF}),
\beglab{ineqcSFBnormF}
\abs{\cS^F}\leq B^{\abs{\norm F}}
\quad \text{for all}\ens F\in NV(\cN),
\edla
where $B= e^{\gamma+\cB(q)}$ (Theorem~\ref{thm:estimq}) or
$B=e^{\gamma+\cB_{v.f.}(\la)}$ (Theorem~\ref{thm:estimvf}),
with numerical constant $\gamma\defeq \sum_{k\geq 1} \big(\frac{1}{k}-\frac{1}{k+1}\big)\log (k+2)$.
%

Taking for granted these results of
Sections~\ref{sec:vanish}--\ref{sec:estim}, we now show how
they allow one to prove the
convergence of the formal linearizing transformation~$h$ when~$a$ is
convergent, together with quantitative estimates.

\subsection{Holomorphy}

Taking for granted~\eqref{eqrestrictedseries}--\eqref{ineqcSFBnormF},
  %
 a simple majorant series argument will yield

\begin{thm}~\label{thm:holomorphic}

  \begin{enumerate}[(i)]
\item Suppose that $g = 
  R_q \circ (\tmop{id}+a) \in \gGnu$ is a local analytic non-resonant 
  diffeomorphism whose spectrum~$q$ satisfies the Bruno condition.
Then the solution~$h$ to the linearization problem is convergent.
%
%
\item Suppose that $V = V_\lambda + a$ is a local analytic non-resonant vector field
  whose spectrum~$\la$ satisfies the Bruno condition.
Then the solution~$h$ to the linearization problem is convergent.
\end{enumerate}
\end{thm}

\begin{proof}
Suppose that we are given a convergent
$a\in\tgPnu$ 
and consider either the holomorphic
diffeomorphism $g = 
R_q \circ (\tmop{id}+a) \in \gGnu$ or the
holomorphic vector field $V = V_\lambda + a$.
Theorems~\ref{thmformaldiffeonu} and~\ref{thmformalfieldnu} give the
solution $h = (h_1,\ldots,h_\nu) \in \tgGnuid$ to the formal linearization problem
in the form 
%
\beglab{eqhirestr}
h_i(z) = 
z_i + \sum_{T\in\cT(\cN)\cap NV(\cN)} \cS^T D_T(a) z_i
\quad \text{for}\ens i=1,\ldots,\nu
\edla
thanks to~\eqref{eqrestrictedseries}, 
  with $\cS^\bullet \defeq S^\bul(q)$ or $\tilde{S}^\bul(\lambda)$.
%
%
As usual, we expand the components of~$a$ as
\[
  a_i(z_1,\ldots,z_\nu) = \sum_{n\in\cN} a_{i,n} z^{n+e_i}
  \quad \text{for $i=1,\ldots,\nu$.} 
\]
We now define $\ov{a}\defeq(\ov{a}_1,\ldots,\ov{a}_\nu)$ by
\[
\ov{a}_i(z_1,\ldots,z_\nu) \defeq \sum_{n\in\cN} \abs{a_{i,n}}
z^{n+e_i}
  \quad \text{for}\ens i=1,\ldots,\nu,
\]
so that $\ov{a}_i$ is a majorant series of~$a_i$ in the following sense:

\begin{itquote}
  A formal series 
$\dst\psi(z)=\sum_{m\in\N^\nu } \psi_m z^m$ 
with real non-negative
coefficients~$\psi_m$ is said to be a majorant series of a formal
series $\dst \phi(z)=\sum_{m\in\N^\nu } \phi_m z^m \in \Cznu$
if $\abs{\phi_m}\leq \psi_m$ for all $m\in\N^\nu$.
We then use the notation $\phi\preceq \psi$.
\end{itquote}
\smallskip

Lemma~\ref{lemDFbetaiF} implies that
\begla
D_T(a)z_i =  \be_{i,T}(a) \, z^{\norm{T}+e_i}
\quad \text{for each}\ens i\in [\nu], \ens T\in \cT(\cN),
\edla
where the coefficient~$\be_{i,T}$ is a polynomial function with rational
\emph{non-negative} coefficients of $\Set{ a_{1,n}, \ldots, a_{\nu,n}
\given n\in\cN }$.
Since $a_i \preceq \ov{a_i}$, this implies that 
\beglab{ineqmajorantDT}
D_T(a)z_i \preceq D_T(\ov{a})z_i.
\edla
Combining this property with~\eqref{ineqcSFBnormF} and~\eqref{eqhirestr}, 
we get
\beglab{eqhiprec}
h_i(z)\preceq z_i + \sum_{T\in\cT(\cN)\cap
  NV(\cN)}B^{\abs{\norm{T}}}D_T(\ov{a})z_i 
\quad \text{for}\ens i\in [\nu].
\edla
Up to vanishing terms that do not alter the result, 
we recognize in the right-hand side of~\eqref{eqhiprec} the action
on~$z_i$ of the tree expansion $\sum K^\bul D_\bul(\ov a)$ with
$K^F\defeq B^{\abs{\norm{F}}}$.
According to Proposition~\ref{prop:geomexpansion}, this tree expansion
is the composition operator~$C_w$,
where $w =(w_1,\ldots,w_\nu)$ is the composition inverse in
$\tgGnuid$ of the tangent-to-identity formal diffeomorphism
$G(z) \defeq z - B\ii \ov a(Bz)$.
%
%
We have found $h_i\preceq w_i$ for each~$i$.
%
Since~$a$ and thus~$\ov{a}$ are convergent, it follows
that~$G$ is a holomorphic tangent-to-identity diffeomorphism, whose
composition inverse~$w$ is also holomorphic.
%
%
Consequently, each component~$w_i$ is a convergent series, and so
is~$h_i$.
%
\end{proof}


\color{black}

\begin{rem}
  Not only does this proof shows that the formal
  series $h_i(z)$ of~\eqref{eqdefhim} have non-trivial domains of
  convergence, but it also shows that their
  expansions~\eqref{eqdefhievalcontr} and~\eqref{eqdefhievalcontrvf}
  as summations over a family of $\cN$-trees are absolutely convergent
  for $\max_i\abs{z_i}$ small enough, as claimed in \cite{EcNonAbord}.
  There is no contradiction with \cite{Eliasson}, which points out that
  (in the more complicated context of the KAM theorem) there is a
  natural way to expand the solution as a sum over a family of trees
  but this sum is \emph{not} absolutely convergent.
  Indeed, the analog of that phenomenon in the Siegel linearization
  problem is that one can give an alternative, more elementary
  representation of the solution as a so-called mould expansion, that
  is a sum over all words on the alphabet~$\cN$, but that sum fails to
  be absolutely convergent, as discussed in \cite{EcNonAbord}.
\end{rem}

\color{black}


\subsection{A lower bound for the radius of convergence}

A little refinement of the previous proof will 
yield a lower bound for the radius of the
polydisc in which the linearizing transformation is holomorphic.

\label{secthmbounds}
\begin{thm}\label{thm:bounds}
Consider $g = R_q \circ (\tmop{id}+a)$ or $V = V_\lambda + a$
with a non-resonant spectrum satisfying the Bruno condition,
\ie the Bruno series $\cB=\cB(q)$ or~$\cB_{v.f.}(\la)$ is finite,
and $a\in\tgPnu$ holomorphic \black{and bounded} in the polydisc
\begla
\D_b^\nu \defeq \Set[\big]{ z=(z_1,\dots,z_\nu)\in \C^\nu \given |z_i|<b
\;\, \text{for}\;\, i=1,\ldots,\nu }.
\edla
%
%
Then 
the solution~$h$ to the linearization problem is holomorphic in
the polydisc~$\D_{b'}^\nu$ with 
\begla
b' \defeq b \frac{e^{-\gamma -\cB}}{\nu(4 M\nu+2)},
%
%
\edla
where $\gamma\defeq \sum_{k\geq 1} \big(\frac{1}{k}-\frac{1}{k+1}\big)\log (k+2)$ and
$M\defeq\frac{1}{b}\max\Set[\big]{ |a_i(z)| \given z\in \D_b^\nu,\; i=1,\ldots,\nu }$.
\end{thm}

%
The proof will require two
lemmas on holomorphic tangent-to-identity diffeomorphisms in $\nu$
complex variables.


\begin{lem}\label{lem:autoholo}
  Let $\nu\geq 1$ and $m >0$. The automorphism $r_m$ of $\tgGnuid$
  defined by
\begin{equation}
  (r_m \varphi)(z)=\frac{1}{m}\varphi(mz)
  \quad \text{for}\ens \varphi \in \tgGnuid
\end{equation}
induces an automorphism of $\gGnuid$. 
More precisely, if~$\varphi$ is convergent in the polydisc~$\D_b^\nu$,
then $r_m \varphi$ is convergent in the polydisc~$\D_{b/m}^\nu$.
\end{lem}

\begin{proof}
  Obvious.
  %
%
\end{proof}

\begin{lem} \label{lem:holinv}
  Let $\nu\geq 1$ and $\al>0$. Define $\varphi_{\al,\nu}=(\varphi_{\al,\nu,1},\dots,\varphi_{\al,\nu,\nu}) \in \tgGnuid$ by
  \[
    \varphi_{\al,\nu,i}(z)=z_i-\al \sum_{k\geq 2} Z^k
    \quad \text{with}\ens Z=z_1+\dots+z_\nu
  \]
for $i=1,\ldots,\nu$.
Then $\varphi_{\al,\nu}$ is holomorphic in
$U_\nu \defeq \lbrace |Z|<1 \rbrace$ and its composition inverse
$\Psi_{\al,\nu}$ is holomorphic in the polydisc
$\D_{\nu\ii\ka(\al \nu)}^\nu$, where
\begin{equation}   \label{eqdefkappa}
  \ka(\be) \defeq 1+2\be-2\sqrt{\be(1+\be)}
  \quad \text{for $\be>0$.}
\end{equation}
Moreover,
\begin{equation}   \label{eq:lowkappa}
\frac{1}{4\be+1}>\ka(\be)> \frac{1}{4\be+2}
  \quad \text{for all $\be>0$.}
\end{equation}
\end{lem}

\begin{proof}
In view of the properties of the geometric series, each $\varphi_{\al,\nu,i}$ is holomorphic in~$U$ and
\[
  \varphi_{\al,\nu,i}(z)=z_i-\al\frac{Z^2}{1-Z}
  \quad \text{for}\ens z\in U. 
\]
In dimension $\nu=1$, with $z=z_1$, the formal inverse $\Psi_{\al,1}$ of $\varphi_{\al,1}$ satisfies
\[
  \Psi_{\al,1}-\al\frac{\Psi_{\al,1}^2}{1-\Psi_{\al,1}}=z,  
\] 
that also reads
\[
(1+\al)\Psi_{\al,1}^2-(1+z)\Psi_{\al,1}+z=0.
\] 
The discriminant of this polynomial is
\[
\Delta=1-2(1+2\al)z+z^2=\big(\ka(\al) -z\big)\big(\ka^*(\al) -z\big)
\]
with~$\ka(\al)$ as in~\eqref{eqdefkappa} 
  %
  %
and $\ka^*(\al) \defeq \frac{1}{\ka(\al)}=1+2\al+2\sqrt{\al(1+\al)}$.
Since $0<\ka(\al)<\ka^*(\al)$, we have
\[
  \Psi_{\al,1}(z)=\frac{1+z-\big(1-2(1+2\al)z+z^2\big)^{1/2}}{2(1+\al)}
  \quad \text{holomorphic for} \ens z\in \D_{\ka(\al)}^1
\]
(using the principal branch of the square root).

In dimension~$\nu$, the composition inverse
$\Psi_{\al,\nu}=(\Psi_{\al,\nu,1},\dots,\Psi_{\al,\nu,\nu})$
of $\varphi_{\al,\nu}$ is easily computed:
since $\sum_{i\in[\nu]}\ph_{\al,\nu,i}(z)= \ph_{\al\nu,1}(Z)$, the
solution to $\Psi_{\al,\nu}\circ\ph_{\al,\nu}=\tmop{id}$ is given by
\[
  \Psi_{\al,\nu,i}(z)=z_i+\frac{1}{\nu}\big(\Psi_{\al\nu,1}(Z)-Z\big)
  \quad \text{for}\ens i=1,\ldots,\nu.
  %
\]
%
The case $\nu=1$ thus ensures that $\Psi_{\al,\nu}$ is holomorphic
in the domain
$\Set[\big]{ (z_1,\ldots,z_\nu)\given Z=z_1+\dots +z_\nu \in
\D_{\ka(\al\nu)}^1 }$,
which contains the polydisc of radius
$\frac{1}{\nu}\ka(\al\nu)$.

 Finally, one gets the estimates~\eqref{eq:lowkappa} as follows.
For $\be>0$, \black{$\ka(\be)+\ka^*(\be)=2(1+2\be)$, thus}
\[
  \ka(\be)-\frac{1}{4\be+2}=\frac{1}{\ka^*(\be)}-\frac{1}{4\be+2}
  =\black{\frac{4\be+2-\ka^*(\be)}{(4\be+2)\ka^*(\be)}}
  %
  %
  =\frac{\ka(\be)}{(4\be+2)\ka^*(\be)}>0
\]
and 
\[
  \ka(\be)-\frac{1}{4\be+1}=
  \black{ \frac{4\be+1-\ka^*(\be)}{(4\be+1)\ka^*(\be)}
  %
  %
  =\frac{\ka(\be)-1}{(4\be+1)\ka^*(\be)}<0.}
\]

\end{proof}

\begin{proof}[Proof of Theorem~\ref{thm:bounds}]
  Suppose that $a$ is holomorphic \black{and bounded} in the polydisc $\D_b^\nu$.
  %
  \black{Replacing~$b$ by $b_1<b$ arbitrarily close to~$b$, without loss of
  generality we can suppose that~$a$ extends continuously to the
  closure of $\D_b^\nu$.}
%
%
The Taylor expansion of its components being written in the
form~\eqref{eqfwithcN}--\eqref{eqfwithcNa}, the Cauchy inequalities
yield 
\[
  \forall i \in [\nu],\ \forall n\in\cN, \quad |a_{i,n}|\leq
  %
  %
  \frac{1}{b^{|n|+1}}\max_{z\in \ov{\D_b^\nu}}\abs{a_i(z)}\leq  \frac{M}{b^{\abs{n}}}
\]
with $M\defeq\frac{1}{b}\max\Set[\big]{ |a_i(z)| \given z\in \D_b^\nu,\; i=1,\ldots,\nu }$.

We now revisit the majorant series argument of the proof of
Theorem~\ref{thm:holomorphic}:
for $i$ in $[\nu]$, we get
\[
  a_i \preceq
  %
  %
  \sum_{n\in\cN,\ n+e_i\in\N^\nu} \frac{M z^{n+e_i}}{b^{|n|}} =
    \sum_{m\in\N^\nu,\ \abs m\ge2} \frac{M z^m}{b^{|m|-1}} \preceq
  \cA(z) \defeq \sum_{k\geq 2} \frac{M Z^k}{b^{k-1}} 
  \quad\text{with}\ens Z=z_1+\dots+z_\nu
\]
(since the multinomial coefficients appearing in~$Z^k$ are bounded
from below by~$1$).
%
Therefore, with $\cA_\nu(z) \defeq (\cA(z),\dots,\cA(z))$, we obtain
$$
D_T(a)z_i \preceq D_T(\cA_\nu)z_i
\quad \text{for all $i\in [\nu]$ and $T\in \cT(\cN)$.}
$$
Combining this property with~\eqref{ineqcSFBnormF} and~\eqref{eqhirestr}, 
we see that the components of the linearization transformation~$h$ satisfy
\beglab{eqhiprecR}
h_i(z)\preceq w_i \defeq z_i + \sum_{T\in\cT(\cN)\cap NV(\cN)}B^{\abs{\norm{T}}}D_T(\cA_\nu)z_i.
\edla
Up to vanishing terms that do not alter the result, 
we recognize in~$w_i$ 
the action on~$z_i$ of the tree expansion
$\sum K^\bul D_\bul(\cA_\nu)$ with
$K^F\defeq B^{\abs{\norm{F}}}$.
According to Proposition~\ref{prop:geomexpansion}, this tree expansion
is the composition operator of the inverse~$w$ of
%
$G(z) \defeq z - B\ii \cA_\nu(Bz)$.
%
%
Using Lemmas~\ref{lem:autoholo} and~\ref{lem:holinv}, we find that
$G = r_{b\ii B}\varphi_{M,\nu}$, 
thus $w = r_{b\ii B}\Psi_{M,\nu}$, which
%
%
is convergent in the polydisc of radius
\[
\frac{b}{B}.\frac{1}{\nu}\ka(M\nu)>\frac{b}{B}\frac{1}{\nu}\frac{1}{4M\nu+2}.
\]
By~\eqref{eqhiprecR}, the components of~$h$ must be convergent in that polydisc too.
%
%
\end{proof}

\section{Vanishing operators}\label{sec:vanish}

In this section, we point out some fundamental properties of the
coarmoulds $D_\bul(a)$ 
%
associated with all possible $a\in \tgPnu$.

\subsection{The set $NV(\cN) \subset \cF(\cN)$}  \label{secNVcN}

 Recall that 
 $$
 \cN \defeq \Set[\big]{ m-e_i \given m \in \N^\nu, \;\; \abs{m} \ge 2, \;\; i
 \in [\nu] }\subset \Z^\nu
 $$
%
%
and, according to Lemma~\ref{lemDFbetaiF}, when writing an arbitrary $\nu$-tuple
$a\in \tgPnu$ as in~\eqref{eqfwithcN}--\eqref{eqfwithcNa},
with coefficients $(a_{i,n})_{(i,n)\in[\nu]\times\cN}$,
each~$D_F$ is a polynomial in these
coefficients, the variables $z_1,\dots,z_{\nu}$ and the derivations
$\pa_1,\dots,\pa_{\nu}$, with non-negative rational coefficients.

\begin{Def}    \label{defunivvanish}
  We say that $D_F$ is universally vanishing if $D_F(a)=0$ for all
  $a\in \tgPnu$.
  %
  We introduce the notation
\[
  V(\cN) \defeq \Set{ F\in \cF(\cN)\given \text{$D_F$ is universally
      vanishing} },
  \qquad NV(\cN) \defeq\cF(\cN) \setminus V(\cN).
\]
\end{Def}

\smallskip

\begin{exa}   \label{exampleNV}
  Let $\nu\ge3$ and $j\in[\nu]$.
  Then there are infinitely many $n\in\cN$ whose $j$th component
  is~$-1$.
  Pick one such~$n$. We then have
  \beglab{eqannulain}
  a\in\tgPnu  \ens\text{and}\ens i\neq j
  \quad\Longrightarrow\quad a_{i,n}=0
  \edla
  by~\eqref{eqfwithcNa}, whence
  $D_{[n]} = a_{j,n} z^{n+e_j} \pa_j$ by Lemma~\ref{lemdefDbulf}(ii),
  making use of bamboo tree notation~\eqref{eqdefbamboon}.
  Now, pick $k\neq j$ and $n'\in\cN$ whose $k$th component is~$-1$ and
  whose $j$th component is~$0$ (there are infinitely many
  possibilities because $\nu\ge3$).
  Lemma~\ref{lemdefDbulf}(ii) yields
  \begla
  D_{[n',n]} = \sum_{1\le i\le\nu} \big( D_{[n]}(a_{i,n'}
  z^{n'+e_i})\big)\pa_i
  = D_{[n]}(a_{k,n'} z^{n'+e_k})\big)\pa_k = 0
  \edla
  (using the analog of~\eqref{eqannulain} for~$n'$ and the fact that
  the $j$th component of $n'+e_k$ is~$0$, hence the corresponding
  monomial is annihilated by~$\pa_j$ and thus by~$D_{[n]}$).
  Hence $D_{[n',n]}$ is universally vanishing, $[n',n]\in V(\cN)$.  
  \end{exa}

  Here are properties of the set $NV(\cN)$ that are
  easily deduced from the recursive definition of the coarmould:

  \begin{prop}   \label{prop:annul1}
    Let $d\in\N^*$, $T_1,\ldots,T_d\in\cT(\cN)$, $F\defeq T_1\cdots
    T_d$ and $n\in\cN$.

\begin{enumerate}[(i)]
\item If $F \in NV(\cN)$ then, for any $j\in[d]$, $T_j\in NV(\cN)$.
\item If $n\att F \in NV(\cN)$ then $F\in NV(\cN)$.
\item We have $\norm{T}\in \cN$ for all $T \in NV(\cN)\cap \cT(\cN)$.
\item If $n\att F \in NV(\cN)$ then $|n|\geq d-1$.
\end{enumerate}

\end{prop}

\begin{proof}
  The two first points are immediate consequences of Lemma~\ref{lemdefDbulf}.
  %
  For the third point, we consider an arbitrary $T \in NV(\cN)\cap
  \cT(\cN)$. Since
\[
D_T=\sum_{1\le i\le\nu} (D_T z_i)\pa_i
\]
is not universally vanishing, there exists 
$i_0\in[\nu]$
such that the term $D_T z_{i_0}$ is non-zero.
By~\eqref{eqDFbeiF} this term can be viewed as
 a monomial $c_{i_0} z^{\norm{T}+e_{i_0}}$ in the
variables~$z$, where $c_{i_0}$ is a non-zero polynomial in the
coefficients $(a_{i,n})_{(i,n)\in[\nu]\times\cN}$. This ensures that
$\norm{T}+e_{i_0}\in \N^{\nu}$, thus $\norm{T}\in \cN$.

For the last point, 
supposing that $T \defeq n\att F\in NV(\cN)$, the formula
\[
D_T=\sum_{1\le i\le\nu} \big(D_F(a_{i,n} z^{n+e_i})\big)\pa_i
\]
shows that there exists~$i_0$ such that $D_F z^{n+e_{i_0}}$ is a non-zero
monomial in the variables~$z$, which can be expressed as a sum over
$i\in[\nu]^d$ by Lemma~\ref{lemdefDbulf}(iii), the term corresponding
to~$i$ being
\begla
\frac{1}{d_1!\cdots d_r!} 
(D_{T_1} z_{i_1}) \cdots (D_{T_d} z_{i_d})
\pa_{i_1} \cdots \pa_{i_d}(z^{n+e_{i_0}}).
\edla
In particular, there must exist $i\in[\nu]^d$ such that
$\pa_{i_1} \cdots \pa_{i_d}(z^{n+e_{i_0}})\neq0$,
and the total degree in~$z$ of that monomial (\ie the
sum of the degrees in each variable $z_i$) is
$|n|+1-d \in \N$, thus $|n|\geq d-1$.
\end{proof}

\subsection{Admissible cuts and composition of operators}

Admissible cuts arise naturally when composing operators
belonging to the family $(D_F(a))_{F\in\cF(\cN)}$.
%
They will play a crucial role in analytic estimates of the armoulds related to
linearization problems through their combinatorial properties.
We recall here some definitions and properties that can be found in
\cite{FM}.

Recall that, according to Section~\ref{seccNFT}, an $\cN$-forest~$F$ can
be represented by an $\cN$-decorated poset and we then denote by~$V_F$
the underlying set of vertices. For an $\cN$-tree~$T$, we denote
by~$\rho_T$ its root, \ie the unique minimal element of $(V_T,\preceq_T)$.

 \begin{Def}
   %
   \color{black}
   An admissible cut of an $\cN$-forest~$F$ is an antichain of
   $(V_F,\preceq_F)$, \ie a subset~$c$ of~$V_F$ in which no two
   elements are comparable.
   %
   %
   We then denote by $P^c (F) \in \cF(\cN)$ the product of the subtrees $\Tree(\sig,F)$ of~$F$ whose
   roots~$\sig$ are in~$c$, and by $R^c (F)\in\cF(\cN)$ the
   remaining part of~$F$, once these subtrees have been removed.
   %
   %
   \color{black}
   We denote the set of all admissible cuts of~$F$ by $adm(F)$.
 \end{Def}

 For example, if $\sigma$ is a vertex of $F=T_1\cdots T_d \in \cF(\cN)$,
 then $c=\{\sigma \}$ is an admissible cut and
 $P^{\{\sigma \}}(F)=\Tree(\sigma,F)$.
 Note also that there are two ``trivial"
 admissible cuts:
 \begin{itemize}
 \item
(empty cut)\;   $\est\in adm(F)$, $P^\est(F)=\emptyset$ and $R^\est(F)=F$;
 \item
(all minimal elements)\;    $m_F\defeq\{ \rho_{T_1},\dots \rho_{T_d}\}\in adm(F)$,
$P^{m_F}(F)=F$ and $R^{m_F}(F)=\emptyset$.
\end{itemize}

Following \cite{FM} (Theorem 5) we have the product rule:

\begin{prop}
  \color{black}
  For any $\cN$-forests $F_1$ and~$F_2$, we have
  \beglab{eqdefprodrule}
  D_{F_1}(a)D_{F_2}(a) = \sum_{F\in\cF(\cN)}k(F_1,F_2,F)D_F(a),
  \edla
  with $k(F_1,F_2,F) \defeq \card \Set{ c\in adm(F)
    \given P^c(F)=F_1,\ R^c(F)=F_2 }$ (note that only finitely many
  terms of the \rhs\ of~\eqref{eqdefprodrule} are nonzero).  
  %
%
\color{black}
\end{prop}


\subsection{Admissible cuts and vanishing operators}

From the product rule~\eqref{eqdefprodrule}, one can deduce new
informations on universally vanishing operators as in \cite[Thm~8]{FM}.

\begin{thm}\label{th:cutvanish}
Given $F\in NV(\cN)$ and $c\in adm(F)$, both $P^c(F)$ and $R^c(F)$ belong to $NV(\cN)$.
\end{thm}

\begin{proof}
%
  %
Suppose that $F\in NV(\cN)$ and $c\in \adm(F)$. 
  The product rule yields
  \begin{align*}
    D_{P^c(F)}D_{R^c(F)} &= \displaystyle \sum_{G\in\cF(\cN)}k(P^c(F),R^c(F),G)D_G \\[.8ex]
    &= \displaystyle k(P^c(F),R^c(F),F) D_F + \sum_{G\not=F} k(P^c(F),R^c(F),F)D_G,
\end{align*}
where $k(P^c(F),R^c(F),F)$ is positive.
In view of~\eqref{eqDFbeiF}, each term of the right-hand
side is a polynomial in the coefficients
$(a_{i,n})_{(i,n)\in[\nu]\times\cN}$, the variables $z_1,\ldots,z_{\nu}$
and the derivations $\pa_1,\dots,\pa_{\nu}$, with non-negative rational
coefficients, and the first one is non-trivial. 
Therefore, as a linear combination with non-negative coefficients,
$D_{P^c(F)}D_{R^c(F)}$ cannot be trivial, and neither factor can.
\end{proof}

This theorem allows one to recover
Proposition~\ref{prop:annul1}(i)
by considering the admissible cut $\{\rho(T_i)\}$,
since $P^{\{\rho(T_i)\}}(T_1\cdots T_d) = T_i$,
as well as Proposition~\ref{prop:annul1}(ii)
by considering the admissible cut~$m_F$ consisting of all minimal elements
of~$F$,
since $P^{m_F}(n\att F) = F$.
 It also allows us to extend the other assertions of Proposition~\ref{prop:annul1}:

\begin{prop}\label{prop:annul2}
  \black{Let $F\in NV(\cN)$.}
  %
%
  %
\begin{enumerate}[(i)]
\item For any $\sigma \in V_F$, $\hat{\sigma}=\norm{\Tree(\sigma,F)}\in\cN$.
\item \black{If~$F$ is an $\cN$-tree, \ie $F=T\in\cT(\cN)\cap NV(\cN)$,
    then,} for any cut $c\in adm(T)$ such that $c\neq\{\rho_T\}$, $\norm{R^c(T)}\in\cN$.
\item For any cut $c\in adm(F)$, 
  \begin{equation}\label{eq:cutdeg}
\abs{\norm{R^c(F)}}\geq \#c-\deg{F}.
\end{equation}
with the convention $\abs{\norm{\emptyset}}=0$.
\end{enumerate}
\end{prop}

\begin{proof} \textit{(i)} 
  For $\sigma \in V_F$, we have already seen that 
  $P^{\{\sig\}}(F)=\Tree(\sigma,F)$.
  Theorem~\ref{th:cutvanish} thus implies
  $\Tree(\sigma,F)\in NV(\cN)\cap \cT(\cN)$, whence
  $\norm{\Tree(\sigma,F)}\in\cN$ by
  Proposition~\ref{prop:annul1}(iii).
  \medskip

  \noindent
\textit{(ii)} Similarly, for an admissible cut $c\neq\{\rho_T\}$ of~$T$, $R^c(T)$ is
an $\cN$-tree and Theorem~\ref{th:cutvanish} implies that it is in
$NV(\cN)\cap \cT(\cN)$, whence $\norm{R^c(T)}\in\cN$ by
  Proposition~\ref{prop:annul1}(iii).
  \medskip

  \noindent
  \textit{(iii)} The last claim is related to Proposition~\ref{prop:annul1}(iv).
  There is nothing to say when $F=\emptyset$, and the case of trivial
  admissible cuts is easy: 
\begin{itemize}
\item When $c=\emptyset$, $R^c(F)=F$, thus $\abs{\norm{R^c(F)}} \geq \#c-\deg{F}$.
  %
%
\item When~$c$ is the set $m_F$ 
  of all minimal elements of~$F$, $R^c(F)=\emptyset$, 
  thus $\abs{\norm{R^c(F)}}= 0 = \#c-\deg{F}$.
  %
%
\end{itemize} 
For the general case, we proceed by induction on the number of vertices $\#F$ of the
$\cN$-forest~$F$, assuming that~$c$ is non-trivial admissible cut of~$F$.
As already mentioned, inequality~(\ref{eq:cutdeg}) holds when $\#F=0$,
we thus suppose that $\#F\geq 1$,
$F = T_1\cdots T_d \in NV(\cN)$, $d=\deg F\ge1$.

We treat separately the cases $d=1$ and $d\ge2$.
%

\smallskip

 -- If $d=1$, then $F=n\att G$ where $n\in \cN$ and $G\in \cF(\cN)$.
  Being non-trivial, $c$ can also be seen as an admissible cut of~$G$, and
$$
R^c(F)=n \att R^c(G).
$$
By Theorem~\ref{th:cutvanish}, $G\in NV(\cN)$, thus the induction
hypothesis yields $\abs{\norm{R^c(G)}}\geq \#c-\deg{G} $ and, by
Proposition~\ref{prop:annul1}(iv), $|n|\geq \deg{G}-1$.
Therefore,
$$\abs{\norm{R^c(F)}}=|n|+\abs{\norm{R^c(G)}}\geq \deg{G}-1+\#c-\deg{G}=\#c-\deg{F}.$$

\smallskip

 --  If $d\geq 2$, 
  then each~$T_j$ is in $NV(\cN)$ by
  Proposition~\ref{prop:annul1}(ii);
  writing $c=c^1\cup \dots\cup c^d$ with $c^j \in adm(T_j)$ 
  and $R^c(F)=R^{c^1}(T_1)\cdots R^{c^d}(T_d)$, we can apply the
  induction hypothesis to each $R^{c^j}(T_j)$. We get
\[\abs{\norm{R^c(F)}}=\sum_{1\le j\le d}\abs{\norm{R^{c^j}(T_j)}}\geq \sum_{1\le j\le d}(\#c^j-1)=\#c-\deg{F}.
\]
%

\smallskip

The proof is thus complete.
\end{proof}

Note that, with reference to~\eqref{eqdefcFp}, Proposition~\ref{prop:annul2}(i) amounts to
\beglab{eqFpincludesNV}
NV(\cN)\subset \cF^+(\cN),
\edla
\black{as announced in Lemma~\ref{propnullDcFplus}.}
Example~\ref{exampleNV} shows that the inclusion is strict if
$\nu\ge3$. Indeed, for this particular $\cN$-tree $[n',n]\in V(\cN)$,
we have
\begla
[n',n] \in \cF^+(\cN)
\ens\Longleftrightarrow\ens
n'+n \in \cN
\ens\Longleftrightarrow\ens
\black{\text{the $k$th component of~$n$ is positive,}}
\edla
whence $V(\cN)\cap\cF^+(\cN)\neq\emptyset$.

\section{Armould estimates}\label{sec:estim}

For a non-resonant diffeomorphism spectrum
$q = (e^{2 i \pi \la_1},\ldots,e^{2 i \pi \la_\nu}) \in (\C^*)^\nu$,
%
the definition~\eqref{eqRuss} of~$\Om(k)$ can be rewritten
\begin{equation}
\Omega(k)= \min \{1\}\cup \Set[\big]{ d(n\cdot \la,\Z) \given  n \in
  \cN \ \text{such that}\ \abs{n}\leq k }
\quad \text{for}\ens k\in\N^*,
\end{equation}
whereas for a non-resonant vector field spectrum $\la\in\C^\nu$,
%
the definition~\eqref{eqBru} of~$\Om_{v.f.}(k)$ can be rewritten
\begin{equation}
\Omega_{v.f.}(k)= \min \{1\}\cup \Set[\big]{ |n\cdot \la| \given  n
  \in \cN \ \text{such that}\ \abs{n} \leq k }
\quad \text{for}\ens k\in\N^*.
\end{equation}
As promised, in this section we will prove the estimates~\eqref{ineqcSFBnormF} for the
armoulds $S^\bul(q)$ and~$\ti S^\bul(\la)$.

\subsection{A counting lemma}
Our proof will rely on the following combinatorial lemma.

\begin{lem}\label{lem:counting}
  Suppose that $F\in NV(\cN)$. Define $W_0(F)=\tilde{W}_0(F)=V_F$ and, for $k\ge1$,
\begin{align*}
    W_k(F) &\defeq \Set[\Big]{ \sigma\in V_F \given d(\la \cdot
                     \norm{\Tree(\sigma,F)},\Z)< \frac{\Omega(k)}{k+2}
             }, \\[.8ex]
    \tilde{W}_k(F) &\defeq \Set[\Big]{ \sigma\in V_F \given |\la
                             \cdot \norm{\Tree(\sigma,F)}|<
                             \frac{\Omega_{v.f.}(k)}{k+2} }.
\end{align*}
Then, for all $k \in \N$, the cardinalities of these sets are bounded
as follows:
\beglab{ineqcardWk}
 %
%
\# W_k(F)\leq \frac{\abs{\norm{F}}}{k+1}, \qquad
 %
%
\# \tilde{W}_k(F)\leq \frac{\abs{\norm{F}}}{k+1}.
\edla
\end{lem}

\begin{proof}
  The proof is exactly the same in both cases, replacing $\Om(k)$ and
  the map $z \mapsto d(z,\Z)$ by $\Om_{v.f.}(k)$ and the map
  $z\mapsto \abs{z}$; the key properties 
\begin{equation}\label{eq:ineg}
\forall z,z'\in \C,\quad d(-z,\Z)=d(z,\Z) \quad\text{and}\quad d(z+z',\Z)\leq d(z,\Z)+ d(z',\Z)
\end{equation}
 obviously hold true for $z\mapsto \abs{z}$ too. 
We thus focus on the proof in the case of $W_k(F)$.

\smallskip

For $k=0$ the result is obvious:  $\#W_0(F)=\# F$ and, since 
all the decorations $n\in \cN$ have $\abs{n}\geq 1$, 
we get $\abs{\norm{F}}\geq \#F$.

\smallskip

We thus suppose $k\ge1$ and proceed by induction on the size of~$F$.
If $\# F=1$, then $F$ as exactly one vertex $\sigma$ decorated by
$n\in\cN$.
Either $W_k(F)=\emptyset$ and the inequality is obvious,
or $W_k(F)=\lbrace \sigma \rbrace$, in which case
\[
d(\la \cdot n,\Z)< \frac{\Omega(k)}{k+2} <\Omega(k)
\]
and, thanks to the definition of $\Omega(k)$, we get $\abs{n}\geq
k+1$, whence
\[
 \frac{\abs{\norm{F}}}{k+1}=\frac{\abs{n}}{k+1}\geq 1=\#W_k(F).
 \]
 We now suppose that the result is true for any $\cN$-forest in $NV(\cN)$ of
 size less than $s\geq 2$ and consider an $\cN$-forest $F\in NV(\cN)$ of size~$s$.
 
If $d\defeq\deg F\geq 2$, then we write $F=T_1\cdots T_d$
with $T_1,\dots,T_d \in \cT(\cN)\cap NV(\cN)$ by
Proposition~\ref{prop:annul1}(i), and since 
 $W_k(F)=W_k(T_1)\cup \dots \cup W_k(T_d)$
and each $\#T_j<s$, the induction hypothesis yields
\[
\#W_k(F)=\#W_k(T_1)+ \cdots + \#W_k(T_d) \leq \frac{\abs{\norm{T_1}}}{k+1}+ \cdots \frac{\abs{\norm{T_d}}}{k+1}=\frac{\abs{\norm{F}}}{k+1}.
\]

It remains to study the case of an $\cN$-forest of degree~$1$, \ie an
$\cN$-tree $F=T$.
We then have $T=n\att G \in NV(\cN)$ with $\#G<s$. 
If $\rho_T \not\in W_k(T)$, then $W_k(T)=W_k(G)$ and $G\in NV(\cN)$ by
Proposition~\ref{prop:annul1}(ii). In this case,
\[
  \#W_k(T)=\#W_k(G) \leq \frac{\abs{\norm{G}}}{k+1} <
  \frac{\abs{\norm{G}}+\abs{n}}{k+1}=\frac{\abs{\norm{T}}}{k+1}.
\]
For the case $W_k(T)=\lbrace \rho_T \rbrace \cup W_k(G)$, we proceed as
follows.
Consider the admissible cut
$c_k=\lbrace \sigma_1,\dots ,\sigma_r \rbrace \in adm(G)\subset adm(T)$
consisting of the minimal elements of $W_k(G)$ (for the partial order
induced by~$\preceq_G$).
%
It is clear that $W_k(G)\subset V_{P^{c_k}(T)}$ and,
since $P^{c_k}(T)\in NV(\cN)$, the induction hypothesis yields
\[
  \#W_k(G) =
  \#W_k(P^{c_k}(T)) \leq \frac{\abs{\norm{P^{c_k}(T)}}}{k+1}.
\] 
The inequality to be proved is equivalent to
\[
  1+\#W_k(P^{c_k}(T)) \leq
  \frac{\abs{\norm{P^{c_k}(T)}}+\abs{\norm{R^{c_k}(T)}}}{k+1}. \]
It is thus sufficient to check that
\beglab{ineqsuffi}
k+1 \leq \abs{\norm{R^{c_k}(T)}}.
\edla
There are two cases:
\begin{itemize}
\item If the cut $c_k$ has at least $k+2$ elements, then Proposition~\ref{prop:annul2}(iii) ensures that
  \[
    \abs{\norm{R^{c_k}(T)}}\geq \#c_k-\deg{T}\geq k+2-1=k+1
  \]
  and we are done.
\item If the cut $c_k=\lbrace \sigma_1,\dots ,\sigma_r \rbrace$ has
  cardinality $r\leq k+1$, we let $U_j\defeq \Tree(\sigma_j,T)$ for $j\in
  [r]$. As $\sig_j\in W_k(T)$ and $\rho_T\in W_k(T)$, we have
  \[
d(\la \cdot \norm{U_j},\Z)<\frac{\Omega(k)}{k+2}, \quad
d(\la \cdot \norm{T},\Z)<\frac{\Omega(k)}{k+2},
\]
so that, using (\ref{eq:ineg}) and the identity
$\norm{R^{c_k}(T)} = \norm{T}-\norm{U_1} -\cdots -\norm{U_r}$,
\begin{align*}
  d(\la \cdot \norm{R^{c_k}(T)},\Z)
 &\leq d(\la \cdot \norm{T},\Z) +d(\la \cdot \norm{U_1},\Z) +\cdots +d(\la \cdot \norm{U_r},\Z) \\[.8ex]
  &< \displaystyle \frac{r+1}{k+2} \Omega(k)\leq \Omega(k)
\end{align*}
whence, since $\norm{R^{c_k}(T)}\in\cN$ (Proposition~\ref{prop:annul2}(ii)),  $\abs{\norm{R^{c_k}(T)}}\geq k+1$.
\end{itemize}
\end{proof}

We are now ready to give estimates for the armoulds $S^\bul(q)$ and
$\tilde{S}^\bul(\la)$. We shall start with the second armould.

\subsection{The case of $\tilde{S}^\bul(\la)$}

\begin{thm}\label{thm:estimvf} If a non-resonant vector field spectrum
  $\la=(\la_1,\ldots,\la_\nu)$ satisfies the Bruno condition
\[
  \cB_{v.f.}(\la) = \sum_{k=1}^{\infty} \left(
    \frac{1}{k}-\frac{1}{k+1}\right) \log\frac{1}{\Omega_{v.f.}(k)} <
  \infty,
\]
then we have 
\beglab{ineqtiSFBnormF}
  \abs{\tilde{S}^F(\la)}\leq B^{\abs{\norm{F}}}
  \quad \text{for all}\ens F\in NV(\cN)
\end{equation}
with
\beglab{eqdefgamma}
  B=e^{\gamma+\cB_{v.f.}(\la)}, \qquad
  \gamma=\sum_{k\geq 1} \left(\frac{1}{k}-\frac{1}{k+1}\right)\log (k+2).
\edla
\end{thm}

We will make use of the easy

\begin{lem}   \label{easylem}
  Under the hypothesis of Theorem~\ref{thm:estimvf}, the sequence
  $\left(\frac{\Omega_{v.f.}(k)}{k+2}\right)_{k\geq 1}$ is positive
  and decreasing, and there exists a sequence $K_n \uparrow\infty$ of
  positive integers such that
  $\dst\frac1{K_n} \ti L_{K_n} \xrightarrow[n\to\infty]{} 0$,
  where
  \beglab{eqdeftiLk}
  \ti L_k \defeq \log \frac{k+2}{\Omega_{v.f.}(k)}
  \quad\text{for each}\ens k\in\N^*.
  \edla
\end{lem}

\begin{proof}[Proof of Lemma~\ref{easylem}]
  The first claim follows from the fact that the sequence $(\Omega_{v.f.}(k))_{k\geq 1}$ is 
  positive and monotonic non-increasing.
  The second claim amounts to $\liminf \frac1{k} \ti L_k=0$,
  a fact that follows from the Bruno condition:
  $\big(\frac1{k} \ti L_k\big)_{k\ge1}$ is a sequence of positive
  numbers that cannot be bounded from below by any $\La>0$ because
  that would imply $\sum \frac1{k(k+1)} \ti L_k =\infty$,
  whereas the latter series is of same nature as the Bruno series,
  since
  \beglab{eqtiLklogOm}
  \tilde{L}_k = \log (k+2) + \log \frac1{\Omega_{v.f.}(k)}
  \edla
  and $\sum\frac1{k(k+1)} \log (k+2) < \infty$.
  \end{proof}

\begin{proof}[Proof of Theorem~\ref{thm:estimvf}]
  Given a non-empty $\cN$-forest $F\in NV(\cN)$, the first claim of
  Lemma~\ref{easylem} entails
  \[
V_F=\tilde{W}_0(F) \supset \tilde{W}_1(F) \supset \cdots \supset \tilde{W}_k(F) \supset \cdots
%
%
\]
and~\eqref{ineqcardWk} implies that $\tilde{W}_k(F)=\emptyset$ for $k\ge \abs{\norm{F}}$.
We use this in the form of a disjoint union $V_F=\dst \bigcup_{k\in\N}
\tilde{W}_k(F)\setminus \tilde{W}_{k+1}(F)$
(in which only finitely many sets are non-empty): \eqref{eqdeftiSFla} implies
\begla
\log\abs{\tilde{S}^F(\lambda)}
%
%
=\sum_{\sig \in V_F} \log\frac{1}{\abs{\lambda \cdot \norm{\Tree(\sigma,F)}}}\\
= \sum_{k\in \N} \, \sum_{\sigma\in \tilde{W}_k(F)\setminus
  \tilde{W}_{k+1}(F)} \log\frac{1}{\abs{\lambda \cdot
    \norm{\Tree(\sigma,F)}}}.
\edla
With the notation $\ti M_k(F)\defeq \#\ti W_k(F)$ for any $k\geq 0$,
there are exactly $\tilde{M}_k(F)-\tilde{M}_{k+1}(F)$ vertices in
$\tilde{W}_k(F)\setminus \tilde{W}_{k+1}(F)$ and, for any such
vertex~$\sig$,
\begin{equation}
\frac{\Omega_{v.f.}(k+1)}{k+3}\leq |\la \cdot \norm{\Tree(\sigma,F)}|< \frac{\Omega_{v.f.}(k)}{k+2}.
\end{equation}
With the notation~\eqref{eqdeftiLk}, this implies that
\begin{equation}\label{eq:crucialineg}
  \sigma\in  \tilde{W}_k(F)\setminus \tilde{W}_{k+1}(F)
  \quad\Longrightarrow\quad 
  \log\frac{1}{\abs{\lambda \cdot \norm{\Tree(\sigma,F)}}} \leq
  %
  %
  \tilde{L}_{k+1},
\end{equation}
thus
\begin{equation}
  \log\abs{\tilde{S}^F(\lambda)}\leq
  \tilde{B}_F(\lambda) \defeq \sum_{k\in \N}(\tilde{M}_k(F)-\tilde{M}_{k+1}(F)) \tilde{L}_{k+1}.
\end{equation}
Since $\ti M_k(F)=0$ for all $k\ge\abs{\norm F}$,
the Abel transform of the series $\tilde{B}_F(\lambda)$ gives
\[
  \tilde{B}_F(\lambda) =
  %
  %
  \tilde{M}_0(F)\tilde{L}_1+\sum_{k= 1}^{K-1}
  \tilde{M}_k(F)(\tilde{L}_{k+1}-\tilde{L}_k)
  \quad \text{for all}\ens K\ge \abs{\norm F}.
\]
%
%
By Lemma~\ref{easylem}, $\tilde{L}_{k+1}-\tilde{L}_k\geq 0$ for $k\geq 1$, thus
Lemma~\ref{lem:counting} yields
\[
  \tilde{B}_F(\lambda)\leq \abs{\norm{F}}\tilde{L}_1 +\sum_{k= 1}^{K-1}
  \frac{\abs{\norm{F}}}{k+1}(\tilde{L}_{k+1}-\tilde{L}_k)
  %
  %
  \quad \text{for all}\ens K\ge \abs{\norm F}.
\]
Using once again the Abel transform, we get
\begin{equation}
  \log\abs{\tilde{S}^F(\lambda)}\leq \ti B_F(\lambda)
  \leq \abs{\norm{F}}
  \bigg( \frac1 K \ti L_K +
    \sum_{k\geq 1} \bigg(\frac{1}{k}-\frac{1}{k+1}\bigg) \tilde{L}_k \bigg),
\end{equation}
where we choose $K$ in the subsequence $(K_n)$ given by Lemma~\ref{easylem}:
in view of~\eqref{eqtiLklogOm}, letting~$n$ tend to~$\infty$, we get
\[
  \log\abs{\tilde{S}^F(\lambda)}\leq 
\abs{\norm{F}} \,
 \underbrace{\sum_{k\geq 1}\left(\frac{1}{k}-\frac{1}{k+1}\right)\log (k+2)}_\gamma
\, + \, \abs{\norm{F}} \cB_{v.f.}(\la),
\]
%
%
whence~\eqref{ineqtiSFBnormF} follows.
\end{proof}

\subsection{The case of $S^\bul(q)$}

The proof is almost the same, except that we need the following lemma:
\begin{lem}\label{lem:maxmodulus}Given $z\in\C$,
\begin{itemize}
\item if $\abs{\Im z}> \frac{1}{2}$, then $\abs{e^{2i\pi z}-1} > \frac{1}{3}$,
 \item if $\abs{\Im z}\leq \frac{1}{2}$, then $\abs{e^{2i\pi z}-1} \geq d(z,\Z)$.
\end{itemize}
\end{lem}

\begin{proof}

If $\Im z > \frac{1}{2}$, then $\abs{e^{2i\pi z}}=e^{-2 \pi \Im z} < e^{-\pi}<1$ so that
\[
\abs{e^{2i\pi z}-1}>1 -e^{-\pi}>\frac{1}{3}.
\]
In the same way, if $\Im z < -\frac{1}{2}$, then $\abs{e^{2i\pi z}}=e^{-2 \pi \Im z} > e^{\pi}>1$ so that
\[
\abs{e^{2i\pi z}-1}>e^{\pi}-1>\frac{1}{3}.
\]
It remains to prove the second part of the lemma. Let us define the square 
$$C=\lbrace z \in \C \ \abs{\Im z} \leq \frac{1}{2} \text{ and } \abs{\Re z} \leq \frac{1}{2} \rbrace.$$
For $z$ in $C$, $d(z,\Z)=\abs{z}$ and, since $z\mapsto e^{2i\pi z}-1$
is $1$-periodic, the inequality is equivalent to
$$
\forall z\in C, \quad \abs{e^{2i\pi z}-1} \geq \abs{z}.
$$
The meromorphic map $g:z\mapsto \frac{z}{e^{2i\pi z}-1}$ is holomorphic in the interior of $C$ and continuous on $C$ so that we can apply the maximum modulus principle:
$$\forall z\in C,\quad \abs{g(z)}\leq \max_{w\in \partial C}\abs{g(w)}.$$
Let us take a look at the four sides of the square $C$:
\begin{itemize}
\item If $w=\pm \frac{1}{2}+iy$ with $-\frac{1}{2}\leq y \leq \frac{1}{2}$,
$$
\abs{g(w)}=\frac{\sqrt{\frac{1}{4}+y^2}}{\abs{-e^{-2\pi y}-1}}\leq \frac{1}{\sqrt{2}}\frac{1}{e^{-2\pi y}+1}\leq \frac{1}{\sqrt{2}}\frac{1}{e^{-2\pi}+1}\leq 1.
$$
\item If $w=x + \eps i\frac{1}{2}$ with $-\frac{1}{2}\leq x \leq \frac{1}{2}$ and $\eps=\pm 1$,
  \begin{multline*}
\abs{g(w)}^2=\frac{x^2+\frac{1}{4}}{\abs{e^{2i\pi x}e^{-\eps
      \pi}-1}^2}=\frac{x^2+\frac{1}{4}}{e^{-2\eps\pi} -2e^{-\eps
    \pi}\cos(2\pi x) +1} \leq \\[.8ex]
 \frac{1}{2} \frac{1}{e^{-2\eps\pi} -2e^{-\eps \pi} +1}=\frac{1}{2(e^{-\eps \pi}-1)^2}
\end{multline*}
and the reader can check that for $\eps=\pm 1$,
$$\frac{1}{2(e^{-\eps \pi}-1)^2}\leq 1.$$
\end{itemize}
On the square $C$ we have $\abs{g(z)}\leq 1$ and this ends the proof.
\end{proof}

We can now state and prove 

\begin{thm}\label{thm:estimq}
  If a non-resonant diffeomorphism spectrum
  $q = (e^{2 i \pi \la_1},\ldots,e^{2 i \pi \la_\nu})$ satisfies the Bruno condition
\[
\cB(q) = \sum_{k=1}^{\infty} \left( \frac{1}{k}-\frac{1}{k+1}\right) \log\frac{1}{\Omega(k)} < \infty,
\]
then we have
\beglab{ineqSFBnormF}
  \abs{{S}^F(q)}\leq B^{\abs{\norm{F}}}
  \quad \text{for all}\ens F\in NV(\cN)
  \edla
  with
\begin{equation}
  B=e^{\gamma+\cB(q)}
  \quad\text{(with $\gamma$ as in~\eqref{eqdefgamma}).}
  %
%
\end{equation}
\end{thm}

\begin{proof}
  We first observe that the sequence $(\Omega(k))_{k\geq 1}$
  is positive, monotonic non-increasing and bounded from above by~$1$,
  so that the sequence $\left(\frac{\Omega(k)}{k+2}\right)_{k\geq 1}$
  is positive, decreasing and bounded from above by~$\frac{1}{3}$.
  
Following the lines of the proof of Theorem~\ref{thm:estimvf},
we get, for any $F\in NV(\cN)$,
\begin{align*}
\log\abs{{S}^F(q)}&=\log \left| \prod_{\sig \in V_F} \frac{1}{q^{ \norm{\Tree(\sigma,F)}}- 1}\right|\\[.8ex]
& = \log \left| \prod_{\sig \in V_F} \frac{1}{e^{ 2i \pi \la \cdot \hat{\sigma}}- 1}\right| \qquad (\hat{\sigma}=\norm{\Tree(\sigma,F)})\\[.8ex]
&=\sum_{\sig \in V_F} \log \frac{1}{\abs{e^{ 2i \pi \la \cdot \hat{\sigma}}- 1}} \\[.8ex]
&= \sum_{k\in \N} \, \sum_{\sigma\in {W}_k(F)\setminus {W}_{k+1}(F)} \log \frac{1}{\abs{e^{ 2i \pi \la \cdot \hat{\sigma}}- 1}}.
\end{align*}
With the notation $M_k(F)\defeq \#W_k(F)$ for any $k\geq 0$, there are exactly
${M}_k(F)-{M}_{k+1}(F)$ vertices in ${W}_k(F)\setminus W_{k+1}(F)$
and, for any such vertex~$\sig$,
\begin{equation}
\frac{\Omega(k+1)}{k+3}\leq d(\la \cdot \hat{\sigma},\Z).
\end{equation}
We now need an analog of the crucial
inequality~(\ref{eq:crucialineg}).
Thanks to Lemma~\ref{lem:maxmodulus},
\begin{itemize}
\item either $\abs{\Im (\la \cdot \hat{\sigma})}\geq \frac{1}{2}$ and
$$
\abs{e^{ 2i \pi \la \cdot \hat{\sigma}}- 1}>\frac{1}{3}\geq \frac{\Omega(k+1)}{k+3},
$$
\item or $\abs{\Im (\la \cdot \hat{\sigma})}\leq \frac{1}{2}$ and 
$$
\abs{e^{ 2i \pi \la \cdot \hat{\sigma}}- 1}\geq d(\la \cdot \hat{\sigma},\Z)\geq \frac{\Omega(k+1)}{k+3}.
$$
\end{itemize}
This means that, with the notation ${L}_k \defeq \log
\frac{k+2}{\Omega_{v.f.}(k)}$, in all cases,
\begin{equation}\label{eq:crucialinegq}
  \sigma\in  {W}_k(F) \setminus W_{k+1}(F)
    \quad\Longrightarrow\quad 
  \log \frac{1}{\abs{e^{ 2i \pi \la \cdot \hat{\sigma}}- 1}}\leq  \log \frac{k+3}{\Omega(k+1)}={L}_{k+1},
\end{equation}
whence
\begin{equation}
\log\abs{{S}^F(q}\leq \sum_{k\in \N}({M}_k(F)-{M}_{k+1}(F)) {L}_k.
\end{equation}
The end of the proof is exactly the same as for Theorem~\ref{thm:estimvf}.
%
%
\end{proof}

\vspace{.5cm}
\appendix


\centerline{\Large\sc Appendix}
\addcontentsline{toc}{part}{\sc Appendix}

\vspace{.1cm}


\section{A closed formula for the coarmould $D_\bul$}
\label{secappcomputDF}


Let us give ourselves $a \in \tgPnu$, with coefficients~$a_{i,n}$ as
in~\eqref{eqfwithcN}--\eqref{eqfwithcNa}.
In this appendix, we give a closed formula for the coarmould $D_\bul=D_\bul(a)$.
%
%
\begin{prop}   \label{propClosedFormDF}
With any map $i \col R \to [\nu]$, where~$R$ is a finite set, and any $m\in\N^\nu$ we associate a non-negative integer
\[
\Ga(i,m) \defeq 
\begin{cases}
  \dfrac{m_1!}{\big(m_1-\di1\big)!} \cdots
  \dfrac{m_\nu!}{\big(m_\nu-\di\nu\big)!} \in \N^*
     & \text{if $\di1\le m_1,\ldots,\di\nu\le m_\nu$,} \\[1ex]
  \qquad\qquad\qquad 0
     & \text{else}
\end{cases}
\]
with the notation $\di1 \defeq \card\big( i\ii(1) \big)$, \ldots, 
$\di\nu \defeq \card\big( i\ii(\nu) \big)$.
Then, for any $F \in \cF(\cN)$, denoting by $R_F \subset V_F$ the set
of all minimal elements of $(V_F,\, \preceq_F)$
and by $S_F^+(\sig)$ the set of all direct successors
of a vertex $\sig \in V_F$,
we have
\beglab{eqDFGaj}
D_F = \frac{z^{\norm{F}}}{\sym(F)} \sum_{j \colb V_F \to [\nu]}
\bigg( \prod_{\sig\in V_F} \Ga\big( j_{|S_F^+(\sig)}, N_F(\sig)+e_{j(\sig)} \big) 
a_{j(\sig),N_F(\sig)} \bigg) \de_{j_{|R_F}}.
\edla
\end{prop}


Here, the notation for $\de_i$ is the same as in
Lemma~\ref{lemDFbetaiF}.
The coefficients $\Ga(i,m)$ can be interpreted as eigenvalues of the
$0$-homogeneous operator~$\de_i$:
\begla
\de_i = z_1^{\di1} \cdots z_\nu^{\di\nu} 
\pa_1^{\di1} \cdots \pa_\nu^{\di\nu}, \qquad
\de_i \, z^m = \Ga(i,m) z^m.
\edla
%
%
The relation between formulas~\eqref{eqDFGaj} and~\eqref{eqDFbeiF} is
\beglab{eqexplicitbeiF}
\be_{i,F} = \frac{1}{\sym(F)} \sum_{\substack{ j \colb V_F \to [\nu]
    \\ \text{such that}\ j_{|R_F} = i }}
\, \prod_{\sig\in V_F} \Ga\big( j_{|S_F^+(\sig)}, N_F(\sig)+e_{j(\sig)} \big) 
a_{j(\sig),N_F(\sig)}.
\edla
Thus we rediscover this way the fact the coefficients~$\be_{i,F}$ are polynomial functions
with rational non-negative coefficients in the variables
$(a_{i,n})_{(i,n)\in[\nu]\times\cN}$.
%

\begin{proof}[Proof of Proposition~\ref{propClosedFormDF}]
Reasoning by induction on the height, it is sufficient to show that, 
\begin{itemize}
\item if~\eqref{eqDFGaj} holds for~$D_F$ with an $\cN$-forest~$F$, then it holds for
$D_{n\att F}$ with whatever $n\in\cN$,
\item if~\eqref{eqDFGaj} holds for $D_{U_1}, \ldots,D_{U_d}$  with
$\cN$-trees $U_1, \ldots,U_d$, then it holds for $D_{U_1\cdots U_d}$.
\end{itemize}
This is easily proved thanks to points~(ii) and~(iii) of Lemma~\ref{lemdefDbulf}.
\end{proof}




\vspace{.7cm}

\noindent
Fr\'ed\'eric Fauvet\\[1ex]
IRMA, University of Strasbourg et CNRS\\
7, rue Descartes\\
67 084 Strasbourg Cedex, France.\\[.5ex]
email:\,{\tt{frederic.fauvet@math.unistra.fr}}


\vspace{.4cm}

\noindent
Fr\'ed\'eric Menous\\[1ex]
Universit\'e Paris-Saclay, CNRS,\\
Laboratoire de math\'ematiques d'Orsay, UMR8628, \\
91405, Orsay, France. \\[.5ex]
email:\,{\tt{frederic.menous@universite-paris-saclay.fr}}


\vspace{.4cm}

\noindent
David Sauzin\\[1ex]
LTE, CNRS, Observatoire de Paris, \\
PSL University, Sorbonne
Universit\'e,
Paris, France
\\[.5ex]
and Capital Normal University, Beijing, China.\\[.5ex]
email:\,{\tt{David.Sauzin@obspm.fr}}

\end{document}